\documentclass[12pt]{amsart}
\usepackage{latexsym, amssymb, amscd, amsfonts}

 \newlength{\baseunit}               
 \newcount{\numlines}                
\newcommand{\C}{\mathbb{C}}

        \newfont{\hollow}{msbm10 scaled\magstep1}
        
        \newfont{\Bfmit}{eufm10 scaled\magstep1}
        \newcommand{\bfmit}[1]{\hbox{\Bfmit {#1}}}

\newcommand{\cH}{{\mathcal H}}
\newcommand{\cK}{{\mathcal K}}

\newcommand{\cU}{{\mathcal U}}

\newcommand{\cW}{{\mathcal W}}

\newcommand{\cE}{{\mathcal E}}
\newcommand{\cO}{{\mathcal O}}
\newcommand{\cL}{{\mathcal L}}

\newcommand{\cF}{{\mathcal F}}

\newcommand{\cQ}{{\mathcal Q}}
\newcommand{\cV}{{\mathcal V}}

\newcommand{\Quot}{\operatorname{Quot}}
\newcommand{\Vol}{\operatorname{Vol}}

\newcommand{\GL}{\operatorname{GL}}
\newcommand{\SL}{\operatorname{SL}}
\newcommand{\diag}{\operatorname{diag}}

\newcommand{\Id}{\operatorname{Id}}

\newcommand{\SU}{\operatorname{SU}}
\newcommand{\su}{\operatorname{su}}

\newcommand{\Det}{\operatorname{Det}}
\newcommand{\Gr}{\operatorname{Gr}}
\newcommand{\Stab}{\operatorname{Stab}}
\newcommand{\rank}{\operatorname{rank}}

\newcommand{\Hom}{\operatorname{Hom}}

\newtheorem{thm}{Theorem}[section]

\newtheorem{lem}[thm]{Lemma}

\newtheorem{prop}[thm]{Proposition}

\newtheorem{conj}[thm]{Conjecture}

\theoremstyle{definition}
\newtheorem{defn}[thm]{Definition}

\newtheorem{rem}[thm]{Remark}           

\theoremstyle{remark}

\newcommand{\lremind}[1]{{}}
\newcommand{\bremind}[1]{{}}

\newcommand{\cut}[1]{}

\begin{document}
\pagestyle{plain} \title{{ \large{Stable Configurations of Linear
Subspaces and Quotient Coherent Sheaves} }
}
\author{Yi Hu}
\address{\newline Department of Mathematics, University of Arizona, Tucson AZ
\newline
Nankai University, Tianjin, China} \email{yhu@math.arizona.edu}


\maketitle


{\parskip=12pt 


\section{Introduction }
In this paper, we begin to study GIT stability of systems of
geometric objects, using the Hilbert-Mumford numerical criterion
and moment map. Here we focus on linear subspaces and quotient
coherent sheaves.

Consider the product $\Pi_{i=1}^m \Gr(k_i, V \otimes W)$ of the
Grassmannians of $k_i$-dimensional subspaces of $V \otimes W$, on
which $\SL(V)$ acts diagonally, where $V$ and $W$ are two fixed
vector spaces over complex numbers. For a set ${\bf \omega}$ of
positive integers, set
$$L_{\bf \omega}=\otimes_{i=1}^{m} \pi_i^*({\cO}_{\Gr(k_i,V \otimes W)}({\omega_i})).$$
$L_{\bf \omega}$ admits a unique $\SL(V)$-linearization.
 Then, using Hilbert-Mumford numerical criterion, we showed that a
 system of linear subspaces $\{K_i \subset V \otimes W \}$, as a point of
$\Pi_{i=1}^m \Gr(k_i, V \otimes W)$,
 is semistable (resp. stable) with respect to
the $\SL(V)$-linearized invertible sheaf $L_{\bf \omega}$ if and
only if, for all nonzero proper subspace $H \subset V$, we have
$$ {1 \over \dim H} \sum_i \omega_i \dim (K_i \cap (H \otimes W))
\le {1  \over \dim V} \sum_i \omega_i \dim K_i$$ (resp. $<$).
This is Theorem \ref{thm:subspace}, which generalizes Mumford's
Proposition 4.3 of \cite{GIT}, where he treated the case  $\Gr(k,
V )^m$, and Dolgachev's Theorem 11.1 of \cite{Dolgachev}, where he
treated the case of subspaces of $V$\footnote{I thank Igor
Dolgachev for informing me of this (after posting the first
version) and two examples of quotients he constructed in \S 11.3.
See the remark after Theorem \ref{thm:subspace}.}. An equivalent
version of the above criterion is given in Theorem
\ref{thm:subspace}' in terms of systems of $a_i$-dimensional
quotients of $V \otimes W$, as points in $\Pi_{i=1}^m \Gr(V
\otimes W, a_i)$. This alternative formulation is necessary for
the later application to quotient coherent sheaves.

To apply moment map, we assume that $\dim W =1$ and consider the
special case of systems of subspaces in $V$. We showed that a
configuration $\{V_i\} \in \prod_i \Gr(k_i, V)$ is polystable if
and only if $\{V_i\}$ can be (uniquely) {\it balanced} with
respect to a Hermitian metric on $V$.
 Here, a Hermitian metric $h$ on $V$ is said to be a {\it balance metric} for the weighted
configuration of vector subspaces $(\{V_i\}, {\bf \omega})$
if the weighted sum of the orthogonal projections
from $V$ onto $V_i$, for $ 1\le i \le m$, is the scalar operator
${\wp}_{\bf \omega}(\{V_i\})={1 \over \dim V} \sum_i \omega_i k_i$. That is
$$\sum_{i=1}^m \omega_i \pi_{V_i} = {\wp}_{\bf \omega}(\{V_i\})\Id_V$$
where $\pi_{V_i}: V \to V_i \hookrightarrow V$ is the orthogonal
projection from $V$ to $V_i$ and $\Id_V$ is the identity map from
$V$ to $V$. In this case, we also say the weighted configuration
$(\{V_i\}, {\bf \omega})$ is {\it balanced}
 with respect to the metric $h$.
We say $(\{V_i\}, {\bf \omega})$ can be (uniquely) balanced if
there is a (unique) $u \in \SU(V) \backslash \SL(V)$ such that
$(\{u  \cdot V_i\}, {\bf \omega})$ is balanced.

When the configuration $\{V_i\}$ is a so-called $m$-filtration,
the existence of a balanced metric was proved by Totaro
\cite{Totaro} where the term {\it good metric} was used instead.
It was also proved in Klyachko's paper \cite{Kly}. Totaro's
motivation is to use {\it good metric} to give an elementary proof
of G. Faltings and G. W\"ustholz's theorem on the stability of
tensor product \cite{Faltings}. Indeed, we have hoped that the
results obtained here may be used to study some problems on
Diophantine approximations.
This is actually one of our original motivations to investigate
the stability of systems of vector subspaces.

Along the way, we generalize the Gelfand-MacPherson correspondence
(\cite{GM}) from configurations of {\it points} to configurations
of {\it linear subspaces}. More precisely, we show that there is a
one-to-one correspondence between the set of $\GL(V)$-orbits on
the product of the Grassmannians $\Pi_{i=1}^m \Gr(k_i, V)$ and the
set of $\GL(k_1) \times \cdots \times \GL(k_m)$-orbits on the
Grassmannian $\Gr(n, {\Bbb C}^{k_1 + \cdots + k_m})$ where $n=\dim
V$ and $\GL(k_1) \times \cdots \times \GL(k_m) \subset \GL(k_1 +
\cdots + k_m)$ acts on coordinate subspaces block-wise. Then,
following the approach of Kapranov (\cite{Kapranov95}), we prove
that there is a natural one-to-one correspondence between the set
of GIT quotients of $\Pi_{i=1}^m \Gr(k_i, V)$ by the diagonal
action of  $\GL(V)$ and the set of GIT quotients of $\Gr(n, {\Bbb
C}^{k_1 + \cdots + k_m})$ by the  action of $\GL(k_1) \times
\cdots \times \GL(k_m)$.  It should follow from here that there is
also an isomorphism between the Chow quotients of the two actions
(cf. Theorem 3.6 of \cite{Hu2003}). When $k_1 = \cdots = k_m =1$,
$\GL(k_1) \times \cdots \times \GL(k_m)$ becomes a maximal torus
of $\GL(k_1 + \cdots + k_m, {\Bbb C})$. And in this case, the
above correspondence becomes the usual Gelfand-MacPherson
correspondence. The case of a product of $\Gr(2,{\Bbb C}^4)$ was
already treated by P. Foth and G. Lozano in \cite{Foth}. After
posting this paper on ArXiv, Ciprian Borcea e-mailed me that his
paper \cite{Borcea} contains a generalization of the
Gelfand-MacPherson correspondence, at birational level,  to flag
configurations.

In addition, by combining the generalized GM correspondence and
the isomorphism between $\Gr(n, {\Bbb C}^{k_1 + \cdots + k_m})$
and $\Gr(k_1 + \cdots + k_m - n, {\Bbb C}^{k_1 + \cdots + k_m})$,
we obtain a {\it generalized Gale transform} from  configurations
of subspaces in $\Pi_{i=1}^m \Gr(k_i, \C^n)$ to configurations of
subspaces in $\Pi_{i=1}^m \Gr(k_i, \C^{k_1 + \cdots + k_m - n})$.
The duality is well-defined up to linear transformations. This
seems to be what was suggested by Eisenbud and Popescu in
\cite{EP}. What is the geometric significance of this duality?
This is a question worth pursuing.

To compute the $\GL(V)$-ample cone of $\Pi_{i=1}^m \Gr(k_i, V)$ or
equivalently the $\GL(k_1) \times \cdots \times \GL(k_m)$-ample
cone of $\Gr(n, {\Bbb C}^{k_1 + \cdots + k_m})$, we  introduce a
new polytope, the {\it diagonal hypersimplex} or {\it
subhypersimplex},  which generalizes the usual hypersimplex (\S
5.2). As an interesting observation, we found that some diagonal
hypersimplexes provide natural examples of $G$-ample cones {\it
without} any top chambers. Not many examples of this sort are
previously known (cf. the Appendix of \cite{DH}).

Finally, as an application, we consider systems of quotient
coherent sheaves.  Let $X$ be a projective scheme (possibly
singular) over the field of complex numbers. Let $\{\cE_i\}$ ($1
\le i \le m$) be a system of (quotient) coherent sheaves  over
$X$, realized as a point in the product of certain Quot schemes
$\Quot(V \otimes \cW,P_i)$ over $X$,  where $V$ is a vector space
and $\cW$ is a coherent sheaf. The group $\SL(V)$ of special
linear transformations acts diagonally on the total product space.
On the product space, there is a $\SL(V)$-linearization
${\cL}_{\bf \omega}$ associated to any given
 set of positive weights  ${\bf \omega} = \{\omega_i\}$
via the Grothendieck embeddings of the corresponding Quot schemes.
We prove that $\{\cE_i\}$ is GIT semistable (resp. stable) with
respect to
 the $\SL(V)$-linearized invertible sheaf ${\cL}_{\bf \omega}$
if and only if for every proper linear subspace $H$ of $V$,
$${1 \over \dim V} \sum_i \omega_i \chi(\cE_i(k))
\le {1 \over \dim H} \sum_i \omega_i \chi(\cF_i(k)) $$ (resp. $<$)
where $\cF_i$ is the subsheaf of $\cE_i$ generated by $H \otimes
\cW$, and $\chi (\bullet)$ is the Euler characteristic. (See \S
\ref{sheaves} for more details.)

Using the relation between GIT stability and the vanishing of
moment map, we proved, in the special case of subbundles of the
trivial bundle $V$, that a configuration $\{\cE_i\}$ of vector
subbundles in $\Pi_{i=1}^m \Quot(V,P_i)$ is polystable if and only
if $\{\cE_i\}$ can be (uniquely) balanced. Here we say that the
configuration $\{\cE_i\}$ of vector subbundles in $\Pi_{i=1}^m
\Quot(V,P_i)$ is balanced if
$$\sum_{i=1}^m \omega_i \int_X A_i(x)A_i^*(x) d V = {\wp}_{\bf \omega}(\{\cE_i\}) \Vol(X) I$$
where $A_i(x)$ is a matrix  representation of $(\cE_i)_x \subset
{\Bbb C}^N$ whose columns form  an orthonormal basis for
$(\cE_i)_x$ ($1 \le i \le m$), $I$ is the identity matrix,   $\Vol
(X)$ is the volume of $X$, and ${\wp}_{\bf \omega} (\{\cE_i\}) =
\sum_i \omega_i {\rank (\cE_i) \over N}$. We say $\{\cE_i\}$ can
be (uniquely) balanced if there is a (unique) element $u \in
\SU(N)\backslash SL(N)$ such that $\{u \cdot \cE_i\}$ is balanced.

When the system consists of a single vector bundle (i.e., $m=1$)
over a smooth projective variety, the above becomes a differential
geometric criterion for the Gieseker-Simpson stability, which is
originally due to  Wang (\cite{Wang}) and Phong-Sturm
(\cite{PhS}). Similar circle of ideas appeared earlier in  the
papers of Zhang (\cite{zhang}) and Luo (\cite{Luo}).

The outcome of this paper relies on the ideas of many other people
in their earlier works, my sole contribution  is to generalize
them to systems of vector subspaces and coherent sheaves, in the
hope that they will be used in future applications and references.
The use of balance metrics was inspired by Totaro \cite{Totaro},
Klyachko \cite{Kly}, and by Wang (\cite{Wang}),
 Phong-Sturm (\cite{PhS}) and the earlier papers
of Zhang (\cite{zhang}) and Luo (\cite{Luo}); The GIT
constructions of the moduli spaces of stable configurations of
coherent sheaves follow very closely the approach of Simpson
(\cite{Simpson94}); The generalized GM correspondence obviously
plainly  follows Gelfand-MacPherson (\cite{GM}) and Kapranov
(\cite{Kapranov95}); The author benefited from the conversations
with P. Foth and W.-P. Li, and from the correspondence with I.
Dolgachev and C. Simpson. I thank them all. Financial support and
hospitality from Harvard University and Professor S.-T. Yau, from
NCTS Taiwan and Professor C.L Wang, and from Hong Kong UST and
Professors W.-P. Li and Y. Ruan are gratefully acknowledged. The
research is partially supported by NSF and NSA. The paper was
finished in early 2003.


\section{Configuration of subspaces and quotients of tensor product}
\label{tensorproduct}

Throughout the paper, we will work over the field of complex
numbers.  Let $V$ and $W$ be two vector spaces. Consider the
product of the Grassmannians
$$\Pi_{i=1}^m \Gr(k_i, V \otimes W).$$
The group $\SL(V)$ acts diagonally on $\Pi_{i=1}^m \Gr(k_i, V
\otimes W)$ by operating on the factor $V$. We will study the GIT
of this action.



\subsection{Stability Criteria}


To proceed we need a lemma.

\begin{lem}
\label{lambda-s} Let ${\bf q}$ be the vector $(q_1, \ldots, q_n)$
such that
$$(\star) \;\;\; q_1 \ge q_2 \ldots \ge q_n, \;\; \hbox{and} \;\;  q_1 + \ldots + q_n = 0.$$
Let ${\bf q}_s $ be the vector $(q_1, \ldots, q_n)$ such that
$$ q_1 = \ldots = q_s = n-s, q_{s+1} = \ldots = q_n = -s, \;\; \hbox{for} \;\; s= 1, \ldots, (n-1).$$
Then ${\bf q}$ is a linear combination of ${\bf q}_s$, $s= 1,
\ldots, (n-1)$, with nonnegative coefficients.
\end{lem}
\proof Indeed, one can check that
$${\bf q} = {{q_1 - q_2} \over {n}} {\bf q}_1 + \ldots + {{q_{n-1} - q_n} \over {n}} {\bf q}_{n-1}.$$
\endproof

Let ${\bf \omega}=\{\omega_1, \ldots, \omega_m \}$ be a set of
positive integers,  and
$$L_{\bf \omega}=\otimes_{i=1}^{m} \pi_i^*({\cO}_{\Gr(k_i, V \otimes W)}({\omega_i}))$$
be the ample line bundle over $\Pi_{i=1}^m \Gr(k_i, V \otimes W)$
associated with ${\bf \omega}$, where $\pi_i$ is the projection
from the product space to the $i$-th factor. This line bundle has
a unique $\SL(V)$-linearization because $\SL(V)$ is semisimple
(\cite{GIT}).

We refer the reader to consult \cite{GIT} for the definition of
GIT stability and for the Hilbert-Mumford numerical criterion.

\begin{thm}
\label{thm:subspace} A system of linear subspaces $\{K_i \subset V
\otimes W \}$ as a point of $\Pi_{i=1}^m \Gr(k_i, V \otimes W)$
 is semistable (resp. stable) with respect to
the $\SL(V)$-linearized invertible sheaf $L_{\bf \omega}$ if and
only if, for all nonzero proper subspace $H \subset V$, we have
$$ {1 \over \dim H} \sum_i \omega_i \dim (K_i \cap (H \otimes W))
\le {1  \over \dim V} \sum_i \omega_i \dim K_i$$ (resp. $<$).
\end{thm}

\proof Choose a basis ${v_1, \ldots, v_n}$ of $V$ such that $H =
\hbox{span}\{v_1, \ldots, v_s\}$. Set $H_i= \hbox{span}\{v_1,
\ldots, v_i\}$. In particular, we have $H_n=V$ and $H_s=H$.

Let ${w_1, \ldots, w_d}$ be a basis for $W$. We list the basis of
$V \otimes W$ made of $v_i \otimes w_j$ as
$$\{v_1\otimes w_1, v_1 \otimes w_2, \ldots, v_n \otimes w_d\}.$$
Let $E_i$ ($1 \le i \le nm$) be spanned by the first $i$ vectors
in the above basis.

Let $K$ be any subspace of $V \otimes W$. Then for any integer $1
\le j \le  k=\dim K$, there are integers
$l_j$ such that
$$\dim K \cap E_{l_j-1} = j-1, \; \dim K \cap E_{l_j} = j.$$

Under the basis $\{v_1\otimes w_1, v_1 \otimes w_2, \ldots, v_n
\otimes w_d\}$, $K$ can be represented by a matrix
$$
\left(
\begin{array}{ccccccc}
            a_{11} & \cdots & a_{1l_1} & 0 & \cdots & 0 & 0 \cdots 0 \\
            a_{21} & \cdots & \cdots & a_{2l_2} & \cdots & 0 & 0 \cdots 0 \\
          \vdots & \vdots & \vdots & \vdots & \vdots & \vdots & \vdots \\
            a_{k1} & \cdots & \cdots & \cdots & \cdots & a_{kl_k} & 0 \cdots 0
   \end{array}
 \right)
$$
In the Plucker embedding, one sees that
$$p_{i_1\cdots i_k} (K) = 0   \;\;\hbox{if}\;\; i_j > l_j $$
$$p_{l_1\cdots l_k} (K) \ne 0 .$$

We now apply the above to all $K_i$ ($1 \le i \le m$) and let
$l_j^{(i)}$ be the numbers associated to $K_i$.

Next, consider the one-parameter subgroup $\lambda (t)$ of
$\SL(V)$ defined by a vector ${\bf q} = (q_1, \ldots, q_n)$ as a
diagonal matrix
$$\lambda (t) = \diag (t^{q_1}, \ldots, t^{q_n})$$
with $$q_1 + \ldots + q_n = 0.$$ By permutation if necessary, we
can further assume that
 $$q_1 \ge q_2 \ldots \ge q_n.$$
Let each $q_i$ repeat $m$ times, we obtain a new diagonal matrix
$$\lambda' (t) = \diag (t^{q_1'}, \ldots, t^{q_{mn}'}).$$
Under this convention and from the matrix representations of $K$,
we see that
$$p_{i_1\cdots i_k} (\lambda (t) K) = t^{q_{i_1}' + \ldots q_{i_k}'} p_{i_1\cdots i_k}(K).$$
Hence by the minimality of the numerical function we obtain
$$\mu^{L_{\bf \omega}}(\{K_i\}, \lambda)
= \sum_{i=1}^m \omega_i \sum_{j=1}^{k_i} q_{l_j^{i}}'.$$

Using the fact that $\dim K_i \cap E_j - \dim K_i \cap E_{j-1}$
equals $0$ when $j \ne l_j^{(i)}$ and equals 1 otherwise, we can
rewrite
$$\mu^{L_{\bf \omega}}(\{K_i\}, \lambda)
= \sum_{i=1}^{m}\omega_i (\sum_{j=1}^{mn} q_{j}'(\dim K_i \cap E_j
- \dim K_i \cap E_{j-1}) ).$$ (Note here that $\mu^{L_{\bf
\omega}}(\{K_i\}, \lambda)$ is linear in $(q_1, \ldots, q_n)$.
This observation will be useful later.)

Now replace  $\lambda$ by the one-parameter subgroup $\lambda_s$
defined by ${\bf q}_s$ ($1 \le s \le (n-1)$, see Lemma
\ref{lambda-s}), then we have
$$\mu^{L_{\bf \omega}}(\{K_i\}, \lambda_s) =
\sum_{i=1}^{m}\omega_i (\sum_{j=1}^{sm} (n-s)(\dim K_i \cap E_j -
\dim V_i \cap E_{j-1}) )$$
$$- \sum_{i=1}^{m}\omega_i (\sum_{j=sm+1}^{nm} s (\dim K_i \cap E_j - \dim V_i \cap E_{j-1}) ).$$
After cancelation, we obtain
$$\mu^{L_{\bf \omega}}(\{K_i\}, \lambda_s) =
\sum_{i=1}^{m} \omega_i ( (n-s)\dim K_i \cap E_{sm} - s(\dim K_i
\cap E_{mn} - \dim K_i \cap E_{sm}) ).$$ That is,
$$\mu^{L_{\bf \omega}}(\{K_i\}, \lambda_s) =
\sum_{i=1}^{m}\omega_i ( n \dim K_i \cap E_{sm} - s\dim K_i \cap
E_{mn}  ),$$ Noting that $E_{sm} = H \otimes W$ and $E_{mn} = V
\otimes W$, we have
$$\mu^{L_{\bf \omega}}(\{K_i\}, \lambda_s) = \dim V \sum_{i=1}^{m}\omega_i \dim K_i \cap (H \otimes W) -
\dim H \sum_{i=1}^{m}\omega_i  \dim K_i.$$

Now if $\{K_i\}$ is $L_{\bf \omega}$-semisimple (resp. simple),
then
$$\mu^{L_{\bf \omega}}(\{K_i\}, \lambda_s) \le 0 $$
(resp. $<$) which is the same as that
$$ {1 \over \dim H} \sum_i \omega_i \dim (K_i \cap (H \otimes W))
\le {1  \over \dim V} \sum_i \omega_i \dim K_i $$ (resp. $<$).

Conversely, if the inequality
$$ {1 \over \dim H} \sum_i \omega_i \dim (K_i \cap (H \otimes W))
\le {1  \over \dim V} \sum_i \omega_i \dim K_i$$ holds for all
$H$, but $\{K_i\}$ is not $L_{\bf \omega}$-semistable. Then there
is one-parameter subgroup $\lambda$ such that $\mu^{L_{\bf
\omega}}(\{K_i\}, \lambda) > 0$. By conjugation and permutation,
we can assume that the vector ${\bf q}$ that defines $\lambda$
satisfies ($\star$) (see Lemma \ref{lambda-s}). Note that such a
vector ${\bf q}$ is a linear combination of ${\bf q}_s$ ($1 \le s
\le n-1$) with non-negative coefficients. Note also from the above
that $\mu^{L_{\bf \omega}}(\{K_i\}, \lambda)$ is linear in ${\bf
q}$. Hence there must exists $s$ ($1 \le s \le n-1$) such that
$\mu^{L_{\bf \omega}}(\{K_i\}, \lambda_s)> 0 $, but this is
equivalent to that
$$ {1 \over \dim H} \sum_i \omega_i \dim (K_i \cap (H \otimes W))
> {1  \over \dim V} \sum_i \omega_i \dim K_i$$
  for some vector subspace $H$, a contradiction.

Similarly, if the strict inequality $$ {1 \over \dim H} \sum_i
\omega_i \dim (K_i \cap (H \otimes W)) < {1  \over \dim V} \sum_i
\omega_i \dim K_i$$ holds for all $H$, then by the above $\{K_i\}$
is $L_{\bf \omega}$-semistable. Assume that it is not stable. Then
there is a $\lambda$ that satisfies ($\star$) of Lemma
\ref{lambda-s} such that $\mu^{L_{\bf \omega}}(\{K_i\},\lambda) =
0$. Then the same arguments as above plus  that we already know
$\mu^{L_{\bf \omega}}(\{K_i\}, \lambda_s) \le 0 $ will yield that
there is $s$ ($1 \le s \le n-1$) such that $\mu^{L_{\bf
\omega}}(\{K_i\}, \lambda_s) = 0 $, but this is equivalent to that
$$ {1 \over \dim H} \sum_i \omega_i \dim (K_i \cap (H \otimes W))
={1  \over \dim V} \sum_i \omega_i \dim K_i$$
  for some vector subspace $H$, a contradiction.

This completes the proof.
\endproof

In the case of systems of linear subspaces of $V$, Dolgachev in
Theorem 11.1, \cite{Dolgachev2003}  already provided a proof of
the stability criterion. More interestingly, \S 11.3 of
\cite{Dolgachev2003} contains two nice explicit examples: 4 lines
in ${\mathbb P}^3$ where the quotient is ${\mathbb P}^2$, and, 6
lines in ${\mathbb P}^3$ where the quotient is a  double cover of
a toric space ramified over an explicitly given hypersurface. It
seems that these are the only explicitly known nontrivial examples
of quotients.

Now let us go back to our setups. The above theorem can also be
equivalently stated in terms of quotients. We will use the
notation $\Gr( V \otimes W, a)$ for the Grassmannian of quotient
linear spaces of $V \otimes W$ of dimension $a$.

Let ${\bf \omega}$ be a set of positive integers and
$$L_{\bf \omega}'=\otimes_{i=1}^{m} \pi_i^*({\cO}_{\Gr(V \otimes W, a_i)}({\omega_i}))$$
be the ample line bundle over $\Pi_{i=1}^m \Gr( V \otimes W, a_i)$
defined by ${\bf \omega}$ where $\pi_i$ is the projection from the
product space to the $i$-th factor.

\noindent Theorem \ref{thm:subspace}'. {\it A configuration $\{V
\otimes W \stackrel{f_i}{\rightarrow} U_i \rightarrow 0\}$ as a
point of $\Pi_{i=1}^m \Gr( V \otimes W, a_i)$ is semistable (resp.
stable) with respect to the $\SL(V)$-linearized invertible sheaf
$L_{\bf \omega}'$ if and only if, for all nonzero proper subspace
$H \subset V$, we have
$${1 \over \dim V} \sum_{i=1}^m \omega_i \dim U_i \le
{1 \over \dim H} \sum_{i=1}^m \omega_i \dim f_i (H \otimes W) $$
(resp. $<$). In particular, $f_i (H \otimes W) \ne 0$ for some
$i$.}

We note that when $m=1$, this is Simpson's Proposition 1.14 of
\cite{Simpson94},  where it is used to construct the moduli space
of coherent sheaves.


For the action of $\SL(V)$ on $\Gr(V \otimes W, a),$ if there is a
GIT quotient, then it will be {\it unique} because there is only
one ample line bundle over $\Gr(V \otimes W, a)$ up to homothety,
and, this line bundle has a unique $\SL(V)$ linearization. There
could be none, for example, this will be the case when $\dim W=1$.
From now on, we assume that a GIT quotient exists and we use
${\bfmit M}$ to denote this unique quotient variety.

Fix a set of positive numbers  ${\bf \omega} = \{\omega_1, \ldots,
\omega_m\}$ and let ${\bfmit M}_{\bf \omega}$ be the quotient
variety  of the locus of the $L_{\bf \omega}'$-semistable
configurations.

\begin{prop}
\label{projection} Fix an integer $1 \le i \le m$. For
sufficiently large $\omega_i$ (relative to other $\omega_j$), we
have
\begin{enumerate}
\item If a configuration $\{V \otimes  W \rightarrow U_j
\rightarrow 0\}$
 is $L_{\bf \omega}'$-semistable, then its $i$-th component
$V \otimes  W \rightarrow U_i \rightarrow 0$ is also semistable.
\item If the $i$-th component of a configuration $\{V \otimes  W
\rightarrow U_j \rightarrow 0\}$ is stable, then the configuration
$\{V \otimes  W \rightarrow U_j \rightarrow 0\}$ is $L_{\bf
\omega}'$-stable.
\item In particular, there is a surjective
projective morphism from ${\bfmit M}_{\bf \omega}$ to ${\bfmit M}$
$$\pi_i: {\bfmit M}_{\bf \omega} \rightarrow {\bfmit M}.$$
\end{enumerate}
Similar statements  hold in terms of  systems of subspaces.
\end{prop}
\proof For any subspace $H \subset V$, the inequality
$${1 \over \dim V} \sum_j \omega_j \dim U_j \le
{1 \over \dim H} \sum_j \omega_j \dim f_j (H \otimes W) \;
(\hbox{resp} <)$$ holds if and only if
$$\dim H   \dim U_i -\dim V \dim f_i (H \otimes W) \le$$
$${1 \over \omega_i} (\dim H \sum_{j\ne i} \omega_j \dim U_j -\dim V \sum_{j \ne i} \omega_j \dim f_j (H \otimes W)) \; (\hbox{resp} <)$$
holds.

Let $R$ be the right hand term of the last inequality. Then we can
choose a sufficiently large $\omega_i$, such that $|R|<1$. Now if
the inequality ``$\le$'' holds, the left hand term
$$\dim H   \dim U_i -\dim V \dim f_i (H \otimes W)$$ must be nonpositive
since it is an integer. This proves (1).

For (2), if $\dim H   \dim U_i -\dim V \dim f_i (H \otimes W)<0$,
then
$$\dim H   \dim U_i -\dim V \dim f_i (H \otimes W) \le -1.$$
Since we have $|R|<1$, we obtain
$$\dim H   \dim U_i -\dim V \dim f_i (H \otimes W) < R$$ which implies
$${1 \over \dim V} \sum_j \omega_j \dim U_j <
{1 \over \dim H} \sum_j \omega_j \dim f_j (H \otimes W).$$

(3). The existence of the morphism $\pi_i: {\bfmit M}_{\bf \omega}
\rightarrow {\bfmit M}$ follows from (1). The surjectivity follows
from (2).
\endproof


\subsection{Harder-Narasimhan and Jordan-H\"older filtrations}

Theorem \ref{thm:subspace} motivates the following definition. Let
${\bf \omega} = (\omega_1, \ldots, \omega_m)$ be a set of positive
numbers,  called the weights. The {\it normalized total weighted
dimension} of $\cK = \{K_i \} \in \prod_i \Gr(k_i, V \otimes W)$
with respect to ${\bf \omega}$ is defined by
$${\wp}_{\bf \omega} ({\cK}) = {1 \over \dim V} \sum_i \omega_i \dim K_i.$$
For any subspace $H$ of $V$, there is an induced subconfiguration
of linear subspaces in $H \otimes W$
$${\cH} = (K_1 \cap (H \otimes W), \ldots, K_m \cap (H \otimes W))$$
whose {\it normalized total weighted dimension}
 with respect to ${\bf \omega}$ is
$${\wp}_{\bf \omega} ({\cH}) = {1 \over \dim H} \sum_i \omega_i \dim K_i \cap (H \otimes W).$$


\begin{defn}
The configuration ${\cK}$ is ${\wp}_{\bf \omega}$-semistable
(resp. stable) with respect to the weights ${\bf \omega}$ if
$$  {\wp}_{\bf \omega} ({\cH}) \le {\wp}_{\bf \omega} ({\cK}) \;\; ( \hbox{resp}. <)$$
for every subspace $H$ of  $V$.
\end{defn}

Then, Theorem \ref{thm:subspace} can be restated as

\begin{thm}
\label{basic} The configuration ${\cK}$ is GIT semistable (stable)
with respect to the linearized line bundle $L_{\bf \omega}$ if and
only if it is ${\wp}_{\bf \omega}$-semistable (stable) with
respected to the weight set ${\bf \omega}$.
\end{thm}

Let $f: V \rightarrow Q$ be a linear map. Then, the induced map $
V \otimes W \rightarrow Q \otimes W$, still denoted by $f$,
induces a configuration $\{ f(K_i) \}$ of linear subspaces in $Q
\otimes W$. A subconfiguration of $\{K_i\}$ is the one induced
from an inclusion map $i: H \hookrightarrow V$.

\begin{lem}
\label{easylemma} Let $\{K_i\}$ be a configuration of linear
subspaces of $V \otimes W$ and  $$0 \rightarrow F \rightarrow V
\rightarrow Q \rightarrow 0$$ be an exact sequence. Let $\cF$ and
$\cQ$ be the inducing configurations. Then
\begin{enumerate}
\item ${\wp}_{\bf \omega} (\cQ) \ge {\wp}_{\bf \omega} ({\cV})
(\hbox{resp}. >)$ if ${\wp}_{\bf \omega} (\cF) \le {\wp}_{\bf
\omega} ({\cV})  (\hbox{resp}. >)$; \item ${\wp}_{\bf \omega}
(\cF) \ge {\wp}_{\bf \omega} (\cQ) (\hbox{resp}. >)$ if
${\wp}_{\bf \omega} (\cF) \ge {\wp}_{\bf \omega} (\cV)
(\hbox{resp}. >)$.
\end{enumerate}
\end{lem}
\proof We prove (2) and leave (1) for the reader.

Let $F_i$ be $K_i \cap (F  \otimes W)$ and $Q_i$ be the image of
$K_i$ under the map $f: V  \otimes W \to Q  \otimes W$ for all
$i$. If ${\wp}_{\bf \omega} (\cF) \ge {\wp}_{\bf \omega} (\cV) $,
then
$$ {1 \over \dim F} \sum_i \omega_i \dim F_i \ge $$
$${1 \over \dim V} \sum_i \omega_i \dim K_i
= {1 \over \dim V} (\sum_i \omega_i \dim F_i + \sum_i \omega_i \dim Q_i).$$
Hence
$$ {\dim F + \dim Q \over \dim F} \sum_i \omega_i \dim F_i \ge
 (\sum_i \omega_i \dim F_i + \sum_i \omega_i \dim Q_i).$$
Therefore
$$ {\dim Q \over \dim F} \sum_i \omega_i \dim F_i \ge   \sum_i \omega_i \dim Q_i.$$
That is, ${\wp}_{\bf \omega} (\cF) \ge {\wp}_{\bf \omega} (\cQ)$.

The strict inequality can be proved similarly.
\endproof

\begin{defn}
For any configuration $\{K_i\}$ of linear subspaces of $V  \otimes
W$, if there is a  filtration
$$0 = V^0 \subset V^1 \subset \cdots \subset V^h =V$$
with the inducing subconfigurations
$$\{0\}=\{K_i^{(0)}\}\subset \{K_i^{(1)}\} \subset \cdots \subset \{K_i^{(h)} \}=\{K_i\}, \;\; 1 \le i \le m$$
where $K_i^{(l)} = K_i \cap (V^l  \otimes W)$ ($1 \le l \le h$)
such that the quotient configuration $\{K_i^{(l)}/K_i^{(l-1)}\}_i$
($1 \le l \le h$) is ${\wp}_{\bf \omega}$-semistable  and   the
normalized total weighted dimension
$${1 \over \dim V^l/V^{l-1}} \sum_i \omega_i \dim({K_i^{(l)}/ K_i^{(l-1)}}), \;\; 1 \le l \le h$$
is strictly decreasing, then the filtration $$0 = V^0 \subset V^1
\subset \cdots \subset V^k =V$$ or rather the filtered
configuration
$$\{0\} =\{K_i^{(0)} \}\subset \{K_i^{(1)}\} \subset \cdots \subset \{K_i^{(h)}\} =\{K_i\},
\;\; 1 \le i \le m$$ will be called a {\it Harder-Narasimhan
filtration} for $\{K_i\}$.
\end{defn}

\begin{prop}
For every configuration $\{K_i\}$ of linear subspaces of $V
\otimes W$, the Harder-Narasimhan filtration exists and is unique.
\end{prop}
\proof

Let $H$ be a subspace of $V$ such that
$${\wp}_{\bf \omega} (\cH) = {1 \over \dim H} \sum_i \omega_i \dim K_i \cap (H \otimes W)$$
is the maximal. If $H = V$, then $\cK$ is ${\wp}_{\bf
\omega}$-semistable, we are done. Otherwise, by maximality, $\cH$
is ${\wp}_{\bf \omega}$-semistable. Now assume $H_1$ is another
linear subspace such that ${\wp}_{\bf \omega} (\cH_1)$ is maximal,
that is, ${\wp}_{\bf \omega} (\cH_1)$ = ${\wp}_{\bf \omega}
(\cH)$. Then $\cH \oplus \cH_1$ is ${\wp}_{\bf \omega}$-semistable
of ${\wp}_{\bf \omega} (\cH \oplus \cH_1)$ = ${\wp}_{\bf \omega}
(\cH)$. Consider the addition map
$$f: H \oplus H_1 \rightarrow V.$$
Since $H \oplus H_1$ is ${\wp}_{\bf \omega}$-semistable, we have
that the normalized weighted dimension of the kernel $\hbox{Ker}
(f)$ is less than or equal to ${\wp}_{\bf \omega} (\cH)$,
therefore the  normalized weighted dimension of the image $\cH +
\cH_1$ is greater than or equal to ${\wp}_{\bf \omega} (\cH)$ (by
Lemma \ref{easylemma} (1)) and hence equal to ${\wp}_{\bf \omega}
(\cH)$  by the maximality of ${\wp}_{\bf \omega} (\cH)$. This
showed that there is a unique subspace $V^1$ such that ${\wp}_{\bf
\omega} (\cK^1)$ is largest, where $\cK^1$ is the induced
configuration from $V^1$. This constitutes the first step of the
filtration
$$0 \subset V^1 \subset V.$$
Next consider $V/V^1$. If $\cK/cK^1$ is semistable, we are done
because Lemma \ref{easylemma} (2) implies that ${\wp}_{\bf \omega}
(\cK^1)
> {\wp}_{\bf \omega} (\cK/\cK^1)$. If $\cK/\cK^1$ is not semistable,
the above procedure can be applied word for word to produce a
unique linear subspace $V^2$ ($V^1 \subset V^2 \subset V$) with
$\cK^2/\cK^1$ semistable. By Lemma \ref{easylemma} (2) again,
${\wp}_{\bf \omega} (\cK^1) > {\wp}_{\bf \omega} (\cK^2/\cK^1)$
because ${\wp}_{\bf \omega} (\cK^1) > {\wp}_{\bf \omega} (\cK^2)$.
Hence by induction, we will obtain the desired filtration.

The uniqueness is clear from the proof.
\endproof

\begin{defn} Assume that $\{K_i\}$ is ${\wp}_{\bf \omega}$-semistable. If
there is a filtration
$$0 = V^0 \subset V^1 \subset \cdots \subset V^k =V$$
with the inducing subconfigurations
$$\{0\}=\{K_i^{(0)}\}\subset \{K_i^{(1)}\} \subset \cdots \subset \{K_i^{(h)} \}=\{K_i\}, \;\; 1 \le i \le m$$
such that the quotient systems $\{K_i^{(l)}/K_i^{(l-1)}\}_i$ are
${\wp}_{\bf \omega}$-stable with the same  normalized total
weighted dimension ${\wp}_{\bf \omega} (\cV)$, then the filtration
is called a Jordan-H\"older filtration.
\end{defn}

\begin{prop} For any ${\wp}_{\bf \omega}$-semistable $\{K_i\}$,
a Jordan-H\"older filtration  exists.
\end{prop}
\proof A construction goes as follows. If $\cK$ is stable, we are
done. Otherwise, let $H$ be a maximal subspace such that
${\wp}_{\bf \omega} (\cH) = {\wp}_{\bf \omega} (\cK)$. Then $\cH$
must also be semistable. Applying Lemma \ref{easylemma} (1) and
(2), one can check that $\cK/\cH$ is ${\wp}_{\bf \omega}$-stable
and ${\wp}_{\bf \omega} (\cK/\cH) = {\wp}_{\bf \omega} (\cK)$.
Repeat the same procedure to $\cH$, we will obtain a desired
filtration.
\endproof

From the proof one see that a Jordan-H\"older filtration always exists but
depends on a choice of maximal subspaces $H$, hence it needs not to be unique.

\subsection{Splitting and Merging}
\label{splittingMerging}

To conclude \S 2, we will make some elementary observations for
the purpose of  future references. 
Let  $${\cK} = (K_1, \ldots, K_m) \in \Gr(k_1, V  \otimes W)
\times \ldots \times \Gr(k_m, V \otimes W)$$ be a configuration of
vector subspaces of $V={\Bbb C}^n$  with weightes ${\bf \omega} =
(\omega_1, \ldots, \omega_m)$. If for every $i$, $\omega_i = s_i +
t_i$
 where $s_i$ and  $t_i$ are nonnegative integers,
then we can split $K_i$ with weight $\omega_i$ into  $K_i$ with
weight $s_i$ and $K_i$ with weight $t_i$. In this way, we obtain a
new configuration $\widetilde{\cK}$ with new weights $\tilde{\bf
\omega}$. We may call such a process splitting or separation.

Conversely, as opposed to splitting, one may consider ``merging''.
That is, for any configuration of vector subspaces
$\widetilde{\cK}$ with weights $\widetilde{\bf \omega}$, if
$\tilde{K}_i =\tilde{K}_j$ for some $i \ne j$, then we can merge
the two as one and count it with new weight $\omega_i=
\tilde{\omega}_i + \tilde{\omega}_j$. This way, we obtain a new
configuration ${\cK}$ with  new weights ${\bf \omega}$. We may
call such a  process  merging.

Clearly in either splitting or merging, we have that
$${\wp}_{\bf \omega} (\cK) = {\wp}_{\tilde{\omega}} (\widetilde{\cK}),$$
and,  it can also  be easily checked that for any subspace $H
\subset K$,
$${\wp}_{\bf \omega} (\cH) = {\wp}_{\tilde{\omega}} (\widetilde{\cH}).$$

For any weights ${\bf \omega}$, if we write every $\omega_i$ as
the sum $1+ \cdots +1$ ($\omega_i$ many), then we obtain a new
weight set ${\Bbb I} = (1, \ldots, 1)$ which we shall call the
trivial weight. Now, let $M_{\bf \omega}$  denote the GIT quotient
 of $X = \Gr(k_1, V \otimes W) \times \ldots \times
\Gr(k_m, V \otimes W)$ defined by the $\SL(V)$-linearized line
bundle $L_{\bf \omega}$. Then it follows that

\begin{prop}
${\cK}$ is semistable (stable) with respect to ${\bf \omega}$ if
and only if $\widetilde{\cK}$ is semistable (stable) with respect
to $\tilde{\omega}$. Consequently, this induces a closed embedding
from the GIT quotient space $M_{\bf \omega}$ to the corresponding
GIT quotient space $M_{\tilde{\omega}}$. In particular, every
$M_{\bf \omega}$ can be embedded in $M_{\Bbb I}$ as a closed
subvariety.
\end{prop}

\section{Balance metrics and Stability}


\subsection{Polystable configurations}

In this section, we will focus on the special case when $\dim
W=1$. So, let $${\cV} = (V_1, \ldots, V_m) \in \Gr(k_1, V ) \times
\ldots \times \Gr(k_m, V)$$ be a configuration of vector subspaces
of $V={\Bbb C}^n$ and ${\bf \omega} = (\omega_1, \ldots,
\omega_m)$ be a set of positive numbers.

\begin{prop}
\label{directsum} If ${\cV}$ is ${\wp}_{\bf \omega}$-semistable
and is a direct sum $\oplus_{i=1}^l \cH_i$ of   a finite number of
subconfigurations, then all $\cH_i$ and $\cV$  have the same
normalized total weighted dimension. In particular, all $\cH_i$
are also semistable.
\end{prop}
\proof
We first prove the case when $\cV = \cH_1 \oplus \cH_2$.
We have
$$0 \rightarrow \cH_1 \rightarrow \cV \rightarrow \cH_2 \rightarrow 0$$
and
$$0 \rightarrow \cH_2 \rightarrow \cV \rightarrow \cH_1 \rightarrow 0.$$
Since $\cV$ is semistable, ${\wp}_{\bf \omega} (\cH_1) \le {\wp}_{\bf \omega} (\cV)$
and ${\wp}_{\bf \omega} (\cH_2) \le {\wp}_{\bf \omega} (\cV)$.
By Lemma \ref{easylemma}, ${\wp}_{\bf \omega} (\cH_2) \ge {\wp}_{\bf \omega} (\cV)$
and ${\wp}_{\bf \omega} (\cH_1) \ge {\wp}_{\bf \omega} (\cV)$. Hence they are all equal.

In general, write $\cV = \cH_i \oplus (\hbox{the rest})$, by the
case $l=2$, ${\wp}_{\bf \omega} (\cH_i) = {\wp}_{\bf \omega}
(\cV)$ for every $i$.
\endproof

\begin{defn}
A semistable configuration ${\cV} = (V_1, \ldots, V_m)$ is called
polystable if it is a  direct sum of  a finite number of stable
subconfigurations of the same normalized total weighted dimension.
\end{defn}

\begin{prop}
\label{polystable}
$\cV$ is polystable if and only if as a point in the product of the Grassmannians
its orbit is closed in the semistable locus.
\end{prop}
\proof
Suppose that $\cV=\{V_i\}$ is polystable and is the direct sum of stable subconfigurations $\{\cH_q\}$
induced from the decomposition $V = \oplus_q H_q$. Let $\cV(t)$ be a curve in $G \cdot \cV$
for $t$ near $t_0$. Let $\cV(0)$ be the limit  of $\cV(t)$ in the semistable locus
at $t_0$. Then $\cV(0)$ is the direct sum
of  $\{\cH_q (0)\}$ where $\cH_q (0)$ is in the closure of $G \cdot \cH_q$. By Proposition \ref{directsum},
$\cH_q (0)$ is semistable.
 Since $\cH_q$ is stable, $\cH_q (0)$ in the orbit $G \cdot \cH_q$. This means there is a linear isomorphism
$l_q$ of $V$ sending $H_q$ to $H_q(0)$ and inducing isomorphisms
between $H_q \cap V_i$ and $H_q(0) \cap V_i$ for all $i$. Since $V
$ is the direct sum $\oplus_q H_q$, one can build a linear
isomorphism $l$ of $V$ from $l_q|_{H_q}$ (for all $q$), sending
$H_q$ to $H_q(0)$ for all $q$ and inducing isomorphisms between
$H_q \cap V_i$ and $H_q(0) \cap V_i$ for all $i$. Hence $\cV(0) =
\{H_q(0)\}_q$ is in the orbit $G \cdot \cV$.
 This shows that the orbit $G \cdot \cV$ is closed in the semistable locus.

Conversely if $\cV$ is a semistable configuration and
the orbit $G \cdot \cV$ is closed in the semistable locus, we need to show that $\cV$ is polystable.
Let $F \subset V$ be a subspace such that
$\{F \cap V_i\}$ constitute the first step in the Jordan-H\"older filtration of $\cV$.
Choose a basis for $F$ and extend it to a basis for $V$.
Then under this basis we can represent each $V_i$ as an ($n \times k_i$) matrix
$$
 \left(
\begin{array}{ccccccc}
            A_i & B_i\\
            0  & D_i
   \end{array}
 \right)
$$
where $A_i$ generates $F \cap V_i$. Let $d$ be the dimension of $F$ and
 $\lambda (t)$ be a one-parameter subgroup of $\GL(V)$ defined by
$$
\lambda (t) = \left(
\begin{array}{ccccccc}
            t I_d & 0\\
            0  & I_{n-d}
   \end{array}
 \right)
$$
Then we have
$$
\lambda (t) V_i = \left(
\begin{array}{ccccccc}
            t A_i & tB_i\\
            0  & D_i
   \end{array}
 \right)
$$
As $t$ tends to zero, this splits the limit as the direct sum
$$F \cap V_i \oplus Q \cap V_i$$
where $Q$ is spanned by the basis element of $V$ that are not in $F$.
By Proposition \ref{directsum}, the configuration
$\{F \cap V_i \oplus Q \cap V_i\}_i$ is semistable.
Since $G \cdot \cV$ is closed in the semistable locus, this shows that
$\{V_i \}$ and $$\{F \cap V_i \oplus Q \cap V_i\}$$ are in the same orbit.
Repeat the Jordan-H\"older process, this will eventually show that
$\cV$ is polystable.
\endproof

\subsection{Balanced metrics and polystable configurations}

\begin{defn}  A Hermitian metric $h$ on $V$ is said to be a balance metric for the weighted
configuration $({\cV}, {\bf \omega})$  of vector subspaces if the
weighted sum of the orthogonal projections from $V$ onto $V_i$ ($
1\le i \le m$) is the scalar operator ${\wp}_{\bf \omega}({\cV})
\Id_V$. That is
$$\sum_{i=1}^m \omega_i \pi_{V_i} = {\wp}_{\bf \omega}({\cV})\Id_V$$
where $\pi_{V_i}: V \to V_i \hookrightarrow V$ is the orthogonal
projection from $V$ to $V_i$ and $\Id_V$ is the identity map from
$V$ to $V$. In this case, we also say that the weighted
configuration $({\cV}, {\bf \omega})$ is balanced with respect to
the metric $h$. We say $({\cV}, {\bf \omega})$ can be (uniquely)
balanced if there is a (unique) $g \in \SU(V) \backslash \SL(V)$
such that $(g\cdot {\cV}, {\bf \omega})$ is balanced.
\end{defn}

\begin{thm}
\label{b=sfinite} A configuration ${\cV} = (V_1, \ldots, V_m)$ is
polystable with respect to a weight set  ${\bf \omega}$ if and
only if there is  a balance metric on $V$ for the configuration.
\end{thm}

\proof First, it is easy to check that under the linearized line
bundle $L_{\bf \omega}$, the moment map
$$\Phi: \Gr(k_1, V ) \times \ldots \times \Gr(k_m, V) \rightarrow \sqrt{-1} \su (V)$$
of the diagonal action of $\SU(V)$   is given by
$$\Phi({\cV}) = \sum_i \omega_i A_iA^*_i - {\wp}_{\bf \omega}({\cV}) I_n$$
where ${\cV} = (V_1, \ldots, V_m) \in \Gr(k_1, V ) \times \ldots
\times \Gr(k_m, V)$ and $A_i$ is a matrix  representation of $V_i$
such that its columns form an orthonormal basis for $V_i$ ($1 \le
i \le m$). (Here using an orthonormal basis $\{e_1, \cdots, e_n\}$
of $V$, we identify $\su(V)$ with $\su(n)$. Also, using the
Killing form, we identify $\su(n)^*$ with $\sqrt{-1}\su(n)$.)

Assume that ${\cV} = (V_1, \ldots, V_m)$ is polystable with
respect to the weighted ${\bf \omega}$. By Proposition
\ref{polystable}, its orbit in the semistable locus is closed.
Hence by, for example, Theorem 2.2.1 (1) of \cite{DH}, there is an
element $g \in \SL(V)$ such that $\Phi (g \cdot {\cV}) =0$. If $g$
is the identity, this means that
$$\sum_i \omega_i A_iA^*_i = {\wp}_{\bf \omega}({\cV}) I_n$$ which is equivalent to
$$\sum_{i=1}^m \omega_i \pi_{V_i} = {\wp}_{\bf \omega}({\cV}) \Id_V$$
because by a direct computation in linear algebra one can verify that the orthogonal projection
$\pi_{V_i}$ can be identified with the matrix
$A_iA^*_i$ under the identification between $V$ and ${\Bbb C}^n$
(using the orthonormal basis $\{e_1, \cdots, e_n\}$).
That is, the standard hermitian metric $h$ is a balance metric on $V$ for the configuration.
Similarly, when $g$ is not the identity,  the
hermitian metric $gh(\bullet, \bullet) = h(g\bullet,g\bullet)$ is a balance metric on $V$ for the configuration.

Conversely, if there is hermitian metric $h'$ such that it is a
balance metric for the configuration ${\cV} = (V_1, \ldots, V_m)$,
then by scaling, we may assume that $h'$ and $h$ have the same
volume form. Hence there is $g \in \SL(V)$ such that $h' = g h$.
This implies that
$$\Phi(g \cdot {\cV}) = 0.$$
Hence (again by, for example, Theorem 2.2.1 (1) of \cite{DH}), the
orbit through $g \cdot {\cV}$ is closed in the semistable locus.
Therefore by Proposition \ref{polystable}, ${\cV}$ is polystable
with respect to  $L_{\bf \omega}$.
\endproof

This theorem was previously known for the so-called m-filtrations with the trivial weights
 ${\Bbb I}=(1, \ldots, 1)$ and was proved by Klyachko (\cite{Kly})
 and Totaro (\cite{Totaro}).

\subsection{Stability of tensor product}

Of special interest is the so-called $m$-filtration.
A filtration $V^\bullet$ is a weakly decreasing configuration of subspaces
$$V=V^0 \supset V^1 \supset \ldots \supset \{0\}.$$
By a $m$-filtration, we mean   a collection $V^\bullet (s)$ of
filtrations of $V$, for $1 \le s \le m$.

In \cite{Faltings} (cf. also \cite{Totaro}), Faltings and
W\"ustholz defined a  stability for $m$-filtration. Their
definition coincides with our  definition when considering the
$m$-filtration as a configuration of vector subspaces with trivial
weights ${\Bbb I} = (1, \ldots, 1)$.

Conversely, if we treat each $V_i$ as a (trivial) filtration $V
\supset V^i \supset \{0\}$ and use splitting process, then any
configuration of vector subspaces with weights ${\bf \omega}$ can
also be considered as a $m$-filtration with trivial weights ${\Bbb
I}$, and again the two stabilities coincides.

Given two $m$-filtrations $V^\bullet (s)$ and $W^\bullet (s)$, $1 \le s \le m$.
We define the tensor product $(V \otimes W)^\bullet (s)$ as
$$(V \otimes W)^l (s) = \sum_{p + q = l} V^p(s) \otimes W^q(s).$$
If  $V^\bullet (s)$ and $W^\bullet (s)$ have attached weights ${\bf \omega}$
and  ${\bf \omega}'$, then we will use the splitting and
merging method to give $\{(V \otimes W)^l (s)\}$ the induced weights  $\tilde{\bf \omega}$.

\begin{prop}
If $V^\bullet (s)$ and $W^\bullet (s)$ are ${\wp}_\omega$-semistable  and
${\wp}_{\omega'}$-semistable, respectively, then $(V \otimes W)^\bullet (s)$
is ${\wp}_{\tilde{\omega}}$-semistable.
\end{prop}

\proof
This proposition marginally generalizes Theorem 1 of \cite{Totaro}. It also follows
from the proof of \cite{Totaro} using the splitting and merging method to relate
weighted filtrations with unweighted (or trivially weighted) filtrations.
\endproof

A different way to prove this may be done via calculating the moment map of
$\Gr(pq, V \otimes W)$ using the moment map of $\Gr(p,V)$ and $\Gr(q,W)$.

\section{Generalized Gelfand-MacPherson correspondence}

\subsection{Correspondence between orbits}

Choosing a basis of $V$, we can identify $V$ with ${\Bbb C}^n$.
Then a $k$-dimensional vector subspace $E  \subset V \cong {\Bbb
C}^n$ can be represented by a full rank matrix $M$ of size $n
\times k$. The group $G = \SL(n, {\Bbb C})$ acts on $M$ from the
left. The group $G_k = \SL(k,{\Bbb C})$ acts on $M$ from the
right. Two such matrices represent the same vector subspace if and
only if they are in the same orbit of $G_k$. Let $U^0_{n,k}$ be
the space of all full rank matrices of size $n \times k$. Then
$\Gr(k,{\Bbb C}^n)$ is the orbit space $U^0_{n,k}/G_k$.

Now assume that $n < k_1 + \ldots + k_m$.
Given a configuration of vector subspaces
$$(V_1, \ldots, V_m) \in \Gr(k_1, n ) \times \ldots \times \Gr(k_m, n),$$
let $(M_1, \ldots, M_m)$ be their corresponding (representative)
matrices. Now, think of
$$M= (M_1, \ldots, M_m)$$ as a matrix of size $n \times (k_1 + \ldots + k_m)$ and let
$U^0_{n,(k_1, \ldots, k_m)}$ be the space of matrices of size $n \times (k_1 + \ldots + k_m)$
such that $M$ and each of its block matrix $M_j$ is of full rank for $1 \le j \le m$.

There are two group actions on $U^0_{n,(k_1, \ldots, k_m)}$: one is the action of the
group $G = \SL(n, {\Bbb C})$ from the left; the other is the action of the product group
$$ G_{k_1, \ldots, k_m}= \hbox{S}(\GL(k_1, {\Bbb C}) \times \ldots \times \GL(k_m, {\Bbb C}))
\subset \SL(k_1 + \cdots + k_m, {\Bbb C})$$ with each
factor acting on the corresponding block from the right. For simplicity, we sometimes
use $G'$ to denote $G_{k_1, \ldots, k_m}$.
Quotienting  $U^0_{n,(k_1, \ldots, k_m)}$
by the group $G_{k_1, \ldots, k_m}$, we obtain
$$X = \Gr(k_1, n) \times \ldots \times \Gr(k_m, n)$$
with the residual group $G = \SL(n, {\Bbb C})$ acting diagonally as usual;
Quotienting $U^0_{n,(k_1, \ldots, k_m)}$
by the group $G= \SL(n, {\Bbb C})$, we obtain
$$Y = \Gr(n, k_1 + \ldots +  k_m)$$ with the residual group
$G'=G_{k_1, \ldots, k_m}$ acting block-wise.

It follows that

\begin{prop}
There is a bijection between $G$-orbits on $X$ and $G'$-orbits on $Y$.
Indeed, there is a homeomorphism between the (non-Hausdorff) orbit spaces
$X/G$ and $Y/G'$.
\end{prop}

When  $k_1 =k_2= \ldots = k_m = 1$, $X$ is $({\Bbb P}^{n-1})^m$
and $G_{1, \ldots,1}$ is a maximal torus of $\SL(m, {\Bbb C})$. In
this case, the proposition is the  Gelfand-MacPherson
correspondence (\cite{GM}).

\proof The correspondence exists
because each set of orbits
 are in one-to-one correspondence with  $G \times G'$-orbits on
$U^0_{n,(k_1, \ldots, k_m)}$.
\endproof

\subsection{Quotients in stages}

From the previous section one naturally expects that the
correspondence between orbits
 will induce a correspondence between the set of GIT quotients of
$$X = \Gr(k_1, n) \times \ldots \times \Gr(k_m, n)$$
by the  group $G = \SL(n, {\Bbb C})$  and the set of GIT quotients
of
$$Y = \Gr(n, k_1 + \ldots +  k_m)$$  by the group
$G'=G_{k_1, \ldots, k_m}$. The detail of this goes as follows.

First, recall  that any $G$-linearized ample line bundle on  $X=
\Gr(k_1, n) \times \ldots \times \Gr(k_m, n)$ must be of the form
$L_{\bf \omega}$ for some weights ${\bf \omega} $. For the
Grassmannian $Y=\Gr(n, k_1 + \ldots + k_m)$, there is only one
line bundle $L={\cO}_Y(1)$ up to homothety. But the character
group of $\GL(k_1, {\Bbb C}) \times \cdots \times \GL(k_m, {\Bbb
C})$ can be identified with ${\Bbb Z}^m$. That is, each set  ${\bf
\omega}$ of positive integers defines  a character
$$\chi_{\bf \omega}: \GL(k_1) \times \cdots \times
\GL(k_m) \rightarrow {\Bbb C}^*.$$ Let $L(\chi_{\bf \omega})$ be
the ample line bundle ${\cO}_Y(1)$ twisted by the character
$\chi_{\bf \omega}$. $L(\chi_{\bf \omega})$ is linearized for
$\GL(k_1, {\Bbb C}) \times \cdots \times \GL(k_m, {\Bbb C})$  and
hence for its subgroup $G' = \hbox{S}(\GL(k_1, {\Bbb C}) \times
\cdots \times \GL(k_m, {\Bbb C}))$.

\begin{thm}
\label{correspondence} There is a one-to-one correspondence
between the set of GIT quotients of $X = \Gr(k_1, n ) \times
\ldots \times \Gr(k_m, n)$ by the  group $G = \SL(n, {\Bbb C})$
and the set of GIT quotients of  $Y = \Gr(n, k_1 + \ldots +  k_m)$
by the group $G'=G_{k_1, \ldots, k_m}$. More precisely, for any
sequence ${\bf \omega}$ of positive integers, we have a natural
isomorphism between $X^{ss}(L_{\bf \omega})/\!/G$ and
$Y^{ss}(L(\chi_{\bf \omega}))/\!/G'$
\end{thm}

When  $k_1 = \ldots k_m = 1$, the theorem was previously
proved by Kapranov using  the standard Gelfand-MacPherson correspondence.
Here we reproduce his proof in the general case.

\proof First, recall that the coordinate ring of $\Gr(k, {\Bbb
C}^n)$ in the Pl\"uker embedding can be identified with the ring
of polynomials $f$ in matrices $M$ of size $n \times k$ such that
$f(M\cdot g) = f(M)$ for all $g \in \SL(k, {\Bbb C})$. In
particular, we have that the section space $\Gamma (\Gr(k, {\Bbb
C}^n), {\cO}_{\Gr(k, {\Bbb C}^n)}(d))$ can be identified with
$$\{f(M) |f(tM) = t^d f(M), f(M \cdot g) = f(M), g \in \SL(k, {\Bbb C})\}$$
for all integers $d >0$.

Now using the group $\GL(n, {\Bbb C})$ in place of $\SL(k, {\Bbb C})$,
the above has an equivalent but more concise expression as follows. Recall that the character
group of $\GL(k, {\Bbb C})$ can be naturally identified with the group of integers ${\Bbb Z}$.
For any integer $d>0$, let $$\chi_d: \GL(k, {\Bbb C}) \rightarrow {\Bbb C}^*$$
be the corresponding character of $\GL(k, {\Bbb C})$. Then we have
$$\Gamma (\Gr(k, {\Bbb C}^n), {\cO}_{\Gr(k, {\Bbb C}^n)}(d))=
\{f |f(M \cdot g) = \chi_d (g) f(M), g \in \GL(k, {\Bbb C})\}.$$ This is
because the two identities:
 $$f(tM) = t^d f(M) \;\;\hbox{and}\;\; f(M \cdot g) = f(M), g \in \SL(k, {\Bbb C})$$
can be combined together in the single
identity $$f(M \cdot g) = \chi_d (g) f(M), g \in \GL(k, {\Bbb C}).$$

From the above and considering the  ring of polynomials in matrices $M$ of size
$n \times (k_1 + \ldots + k_m)$
one checks that
$$\Gamma(X, L_{\bf \omega}^d) = \{ f(M) | f(M \cdot g)
= \chi_{d{\bf \omega}} (g) f(M), g \in \GL(k_1) \times \cdots \times \GL(k_m)\}$$
and
$$\Gamma(Y, L^d(\chi_{\bf \omega})) = \{ f(M) | f(g' \cdot M) = \chi_{d} (g') f(M), g' \in \GL(n,{\Bbb C})\}.$$

Therefore by taking the projective spectrum of
the invariants of $$A= \oplus_d \Gamma(X, L_{\bf \omega}^d)$$ under the action of the group
$\GL(n, {\Bbb C})$
and by taking the projective spectrum of
the invariants of $$B= \oplus_d \Gamma(Y, L^d(\chi_{\bf \omega}))$$
under the action of the group $$\GL(k_1) \times \cdots \times \GL(k_m),$$
we see that the both quotients
$$X^{ss}(L_{\bf \omega})/\!/G \;\; \hbox{and} \;\; Y^{ss}(L(\chi_{\bf \omega}))/\!/G'$$
can be naturally identified with the projective spectrum of the ring
$$R = \oplus_d R_d$$ where
$$R_d = \{f(M) | f(M \cdot g) = \chi_{d{\bf \omega}} (g) f(M), f(g' \cdot M) = f(M) \}$$
for all $g \in \GL(k_1) \times \cdots \times \GL(k_m), g'\in \SL(n,{\Bbb C})$.
(Note that here we take $g'\in \SL(n,{\Bbb C})$ instead of $g'\in \GL(n,{\Bbb C})$.
This is because the effect of the central part of $\GL(n,{\Bbb C})$ is
already reflected by the scalar matrices in $\GL(k_1) \times \cdots \times \GL(k_m)$.)

This has established the desired correspondence.
\endproof

\section{The Cone of effective linearizations}

As the stability depends on ${\bf \omega}$, so does the moduli. In
this section, we study the $G$-ample cone to pave a way for the
study of the variation of the moduli. In particular, we will
introduce a family of new polytopes: {\it diagonal
hypersimpleces}.

\subsection{Effective linearizations}

Given a linearized line bundle $L$ over $X$, it is called
$G$-effective if $X^{ss}(L) \ne \emptyset$. Not all of $L_{\bf
\omega}$ are $G$-effective. The following should characterize the
effective ample ones.

We will always assume that the group $G= \SL(V)$ acts freely on
generic configuration of linear subspaces, that is, $G$ acts freely on an open
subset of generic points in $\Pi_{i=1}^m \Gr(k_i, V)$.
This should be true when $n < k_1+\cdots +k_m$ and $n^2 \le \sum_i k_i(n-k_i)$.

\begin{conj}
\label{mainconj} Under the above (and perhaps some additional
natural) conditions, we have
\begin{enumerate}
\item $X^{ss}(L_{\bf \omega}) \ne \emptyset$ if and only if $ \omega_i \le {1 \over n} \sum_i k_i \omega_i$ for all $1 \le i \le m$
if and only if $\hbox{max} \{\omega_i \}_i \le {1 \over n} \sum_i k_i \omega_i$;
\item $X^{s}(L_{\bf \omega}) \ne \emptyset$ if and only if $ \omega_i < {1 \over n} \sum_i k_i \omega_i$ for all $1 \le i \le m$
if and only if $\hbox{max} \{\omega_i \}_i < {1 \over n} \sum_i k_i \omega_i$.
\end{enumerate}
The necessary parts of both (1) and (2) are true.
\end{conj}
\proof
 (1). The necessary direction is easy.  Assume that $X^{ss}(L_{\bf \omega}) \ne \emptyset$ and let
$\cV = \{V_i\} \in X^{ss}(L_{\bf \omega})$. We have that for all $W \subset V$,
$${1 \over \dim W} \sum_j \omega_j \dim (V_j \cap W) \le {1 \over n} \sum_i k_i \omega_i.$$
Now for any given $i$, take $W= V_i$, then we obtain
$$\omega_i \le {1 \over \dim W} \sum_j \omega_j \dim (V_j \cap W) \le {1 \over n} \sum_i k_i \omega_i$$
for all $i$.

The necessary part of (2) can be proved similarly.
\endproof

Equivalently,

\begin{conj}
\begin{enumerate}
\item $Y^{ss}(L({\chi_{\bf \omega}})) \ne \emptyset$ if and only if
$ \omega_i \le {1 \over n} \sum_i k_i \omega_i$ for all $1 \le i \le m$
if and only if $\hbox{max} \{\omega_i \}_i \le {1 \over n} \sum_i k_i \omega_i$;
\item $Y^{s}(L(\chi_{{\bf \omega}})) \ne \emptyset$ if and only if $ \omega_i < {1 \over n} \sum_i k_i \omega_i$ for all $1 \le i \le m$
if and only if $\hbox{max} \{\omega_i \}_i < {1 \over n} \sum_i k_i \omega_i$.
\end{enumerate}
The necessary parts of both (1) and (2) are true.
\end{conj}

\subsection{Diagonal  hypersimplex and $G$-ample cone}

The previous conjectures lead to the  discovery of the following
polytope. Setting $x_i =  n \omega_i / \sum_i k_i \omega_i$, then
$x_i$ satisfy $0 \le x_i \le 1$ and $ \sum_i k_i x_i = n.$ Hence
we introduce the polytope
$$\Delta^m_{n,\{k_i\}}=\{(x_1, \ldots, x_m) | 0\le x_i \le 1,\sum_i k_i x_i = n.\}.$$
Recall the standard hypersimplex $\Delta^m_n$ is defined as
$$\Delta^m_n=\{(x_1, \ldots, x_m) | 0\le x_i \le 1,\sum_i  x_i = n.\}.$$
Thus $\Delta^m_{n,\{k_i\}}$ is a subpolytope of $\Delta^{k_1+\cdots +k_m}_n$. In fact,
let $D_{k_1,\ldots,k_m}$ be the diagonal subspace of ${\Bbb R}^{k_1+\cdots +k_m}$ such that the first $k_1$ coordinates
coincide, the next $k_2$ coordinate coincide, and so on, then
$$\Delta^m_{n,\{k_i\}} = \Delta^{k_1+\cdots +k_m}_n \cap D_{k_1,\ldots,k_m}.$$
Clearly $\Delta^m_{n,\{k_i\}}$  is  the hypersimplex $\Delta^m_n$
when all $k_i$ are equal to 1. Hence, it seems reasonable to call
$\Delta^m_{n,\{k_i\}}$ a {\it diagonal hypersimplex} or simply a
{\it generalized hypersimplex}.

Let $G = \SL(V)$ and $G' = S(\GL(k_1, {\Bbb C}) \times \cdots \times \GL(k_m, {\Bbb C}))$.
Let also $C^G(X)$ and $C^{G'}(Y)$ be the $G$-ample cone of $X$ and
$G'$-ample cone of $Y$, respectively. (For the definition and properties of
a general $G$-ample cone, see Definition 3.2.1 and \S3 of \cite{DH}.)

Then, as a corollary of either of the above conjectures, we have

\begin{conj}
\label{gamplecone}
Both $C^G(X)$ and $C^{G'}(Y)$ can be naturally identified with the positive cone over
the generalized hypersimplex $\Delta^m_{n,\{k_i\}}$.
\end{conj}

\subsection{Walls and Chambers}

In \S3 of \cite{DH}, a natural {\it wall and chamber structure} in $C^G(X)$
 is introduced. However it can happen that there are no (top) chambers at all
in $C^G(X)$. Not many examples of this type are previously  known.
Here we produce an interesting one.

Consider the product of $m$-copies of $\Gr(2, {\Bbb C}^4)$,
$$X = \Pi_{i=1}^m \Gr(2, {\Bbb C}^4).$$

\begin{prop} All the above conjectures are true for $\Pi_{i=1}^m \Gr(2, {\Bbb C}^4).$
\end{prop}
\proof We only need to prove it for Conjecture \ref{mainconj}, the
rest follow from this. Take any configuration $\{V_i\} \in
\Pi_{i=1}^m \Gr(2, {\Bbb C}^4)$ such that $V_i \cap V_j = \{0\}$
($i \ne j$). We will check that $\{V_i\}$ is $\wp_{\bf
\omega}$-semistable.  First note that ${ 1 \over \dim V} \sum_i
\omega_i \dim V_i ={1 \over 2} \sum_i \omega_i$. Let $F$ be an
arbitrary proper subspace of $V$. We examine it case by case.

$\dim F =1$. $F$ can intersect non-trivially (i.e., be contained in) only one $V_i$.
Hence we have
$${ 1 \over \dim F} \sum_i \omega_i \dim (F \cap V_i) \le \omega_i \le
{1 \over 2} \sum_i \omega_i.$$

$\dim F =2$. If $\dim F \cap V_i =2$, then $F= V_i$,
hence $${ 1 \over \dim F} \sum_i \omega_i \dim (F \cap V_i) = \omega_i \le
{1 \over 2} \sum_i \omega_i.$$
Otherwise, $\dim F \cap V_i \le 1$ for all $i$, hence
$${ 1 \over \dim F} \sum_i \omega_i \dim (F \cap V_i) \le {1 \over 2} \sum_i \omega_i.$$

$\dim F =3$. If $\dim F \cap V_i \le 1$ for all $i$, then the stability condition is trivially
true. Otherwise, $\dim F \cap V_i =2$ can only be true for only one $i$. In this case,
$${ 1 \over \dim F} \sum_i \omega_i \dim (F \cap V_i) \le {1 \over 3} (2 \omega_i
+ \sum_{j \ne i} \omega_j) $$
$$= {1 \over 3} ( \omega_i + \sum_j  \omega_j) \le {1 \over 2} \sum_j \omega_j.$$
\endproof

An equivalent version of the following proposition already
appeared in Foth-Lozano's paper \cite{Foth} in terms of polygons.

\begin{prop} (\cite{Foth})
For every weight set ${\bf \omega} \in \Delta^m_{4,\{2, \cdots,
2\}}$,
$$X^{ss}(L_{\bf \omega}) \setminus X^s(L_{\bf \omega}) \ne \emptyset.$$ In particular,
there is not any (top) chamber in the $G$-ample cone.
\end{prop}
\proof
Let $F$ be a 2-dimensional subspace. Take a configuration $\{V_i\} \in \Pi_{i=1}^m \Gr(2, {\Bbb C}^4)$
such that $V_i \cap V_j = \{0\}$ ($i \ne j$) and $\dim V_i \cap F =1$ for all $i$.
Then $\{V_i\} $ is semistable for all ${\bf \omega} \in \Delta^m_{4,\{2, \cdots, 2\}}$
by the proof of the previous proposition.
Since
$${ 1 \over \dim F} \sum_i \omega_i \dim (F \cap V_i) = {1 \over 2} \sum_i \omega_i
= { 1 \over \dim V} \sum_i \omega_i \dim V_i,$$ $\{V_i\} \in
X^{ss}(L_{\bf \omega}) \setminus X^s(L_{\bf \omega}) $ for all
$${\bf \omega} \in \Delta^m_{4,\{2, \cdots, 2\}}.$$
\endproof

\begin{rem}
Finally, note that
$$\Delta^m_{4,\{2, \cdots, 2\}}=\{(x_1, \ldots, x_m) | 0 \le x_i \le 1, \sum_{i=1}^m 2 x_i =4\}$$
$$=\{(x_1, \ldots, x_m) | 0 \le x_i \le 1, \sum_{i=1}^m  x_i =2\}$$
which is just the standard hypersimplex $\Delta_2^m$. Recall we
just showed that $\Delta^m_{4,\{2, \cdots, 2\}}$ has no chambers
as $\SL(4,{\Bbb C})$-ample cone of $\Gr(2,{\Bbb C}^4)^m$. However,
$\Delta_2^m$, as $\SL(2,{\Bbb C})$-ample cone of $({\Bbb P}^1)^m$
has natural wall and chamber structure. It would be interesting to
investigate in detail the implications of the above on the problem
of variation of GIT quotients by the two distinct, yet related
actions. Likewise, one should also study the implication of the
identity $$\Delta^m_{kn, \{n, \ldots, n\}} = \Delta^m_k.$$
\end{rem}

\section{Stable Configuration of Coherent Sheaves}
\label{sheaves}

\subsection{Quot scheme and Grothendieck embedding}
Let $X$ be a projective scheme over ${\Bbb C}$ (possibly singular)
with a very ample invertible sheaf $\cO(1)$.
The Hilbert polynomial $p(\cE,k) = \chi (\cE(k))$ is uniquely defined by
the condition that $$p(\cE,k) = \dim H^0(X, \cE(k)),\;\; \hbox{for}\;\; k >>1.$$
Let $d=d(\cE)$ denote the dimension of the support of $\cE$. It is equal to the degree
of $p(\cE,k)$. So,
$$p(\cE,k) = {r \over d!}k^d + {a \over (d-1)!} k^{d-1} + \cdots.$$
Here $r$ is the rank of $\cE$ and $a/r$ is defined to be the slope of $\cE$.
We say $\cE$ is of pure dimension if for any $0 \ne \cF \subset \cE$, we have
$d(\cF) = d(\cE)$.

Fix a vector space $V$ and a coherent sheaf $\cW$ over $X$. Also
fix a (Hilbert) polynomial $P$. We will consider the Quot scheme
$$\Quot(V \otimes  \cW ,P),$$ parameterizing the coherent quotient
sheaves
$$V \otimes  \cW \rightarrow \cE \rightarrow 0$$
such that $p(\cE,k) = P(k)$.

For $k >> 1$, Grothendiek proves that there is an explicit
embedding $\Quot(V \otimes  \cW ,P) \rightarrow  \Gr(V \otimes W,
P(k))$ where $W = H^0(\cW(k))$. Indeed, let ${\cU}$ be the
universal quotient sheaf over $$\Quot(V \otimes \cW, P) \times
X,$$ and $\cL (k)= \Det(p_*({\cU} \otimes q^* \cO_X(k))$ be the
determinant line bundle over $\Quot(V \otimes \cW, P)$ where $p$
and $q$ are the natural projections
\begin{equation*}
\begin{CD}
 \Quot(V \otimes \cW, P) \times X @>{q}>> X  \\
@V{p}VV  \\
 \Quot(V \otimes \cW, P)
\end{CD}
\end{equation*}
Then this is very ample for $k >> 1$ and is the same as the ample line bundle
induced from the embedding into the Grassmaniann (see 1.32 of \cite{Viehweg95}).



\subsection{Stability of configurations of coherent sheaves}

Consider a configuration of coherent quotient sheaves
$$\{V \otimes  \cW \rightarrow \cE_i \rightarrow 0\}_i$$
with $p(\cE_i,k) = P_i(k)$ where $P_i$ are some fixed Hilbert
polynomials. Let $L_{k,i}$ be the linearized ample line bundle on
$\Quot(V \otimes  \cW ,P_i)$ induced from the embedding $$\Quot(V
\otimes  \cW ,P_i) \rightarrow  \Gr(V \otimes W, P_i(k))$$ for
sufficiently large $k$ (we choose $k$ so large that it works for
all $i$). For a given  set of positive integers ${\bf \omega} =
\{\omega_1, \ldots, \omega_m\}$, let $L_{k, {\bf \omega}}$ be the
linearization on
$$\Pi_i \Quot(V \otimes  \cW ,P_i) \subset \Pi_i \Gr(V \otimes W, P_i(k))$$
defined by
$$L_{k, {\bf \omega}}= \otimes_i L_{k,i}^{\omega_i}.$$

We need a simple lemma

\begin{lem} (Theorem 1.19, \cite{GIT}).
\label{inclusionstable}
Let $i: Z \rightarrow \widetilde{Z}$ be a $G$-invariant closed embedding
from a scheme $Z$ to a scheme $\widetilde{Z}$
and $L$ a linearized ample
line bundle over $\widetilde{Z}$. Then
$Z^{ss}(i^*L) = i^{-1} (\widetilde{Z}^{ss}(L))$ and $Z^s = i^{-1}(\widetilde{Z}^s)$.
\end{lem}

This lemma when applied to the Grothendiek embedding will allow us to work directly on
the Grassmannian instead of the Quot scheme.

\begin{thm}
There is an integer $M$ such that for $k\ge M$, the following
holds. Suppose that $\{V \otimes  \cW \stackrel{f_i}{\rightarrow}
\cE_i \rightarrow 0\}$ is a point in
$$\Pi_i \Quot(V \otimes  \cW ,P_i),$$
and for any subspace $H \subset V$, let $\cF_i$ denote the subsheaf of $\cE_i$ generated by
$H \otimes \cW$. Then $\{\cE_i\}$ is semistable (resp. stable) with respect to
the $\SL(V)$-linearization $L_{k,{\bf \omega}}$ if and only if
$${1 \over \dim V} \sum_i \omega_i \chi(\cE_i(k)) \le {1 \over \dim H} \sum_i \omega_i \chi(\cF_i(k))$$
(resp. $<$). In particular, $\chi(\cF_i(k)) >0$ for some $i$.
\end{thm}
\proof
For $k >> 1$, we have the product of the Grothendiek embeddings
$$\Pi_i \Quot(V \otimes  \cW ,P_i) \rightarrow \Pi_i \Gr(V \otimes W, P_i(k))$$
where $W = H^0(\cW(k))$.
Consider the sequences
$$\{ H \otimes \cW  \stackrel{f_i}{\rightarrow} \cF_i \rightarrow 0\}.$$
Let $\cK_i$ be the kernal of $f_i$.
Since all such $H$ runs over a bounded family, so does $\cF_i$. Hence $\cK_i$ also runs
over a bounded family. In particular we may choose $M$ large enough so that
when $k \ge M$, $\chi(\cF_i(k)) = h^0(\cF_i(k))$ and $h^1(\cK_i(k)) =0$ for all such
$\cF_i$ and $\cK_i$. Twist the exact sequence
$$0 \rightarrow \cK_i \rightarrow H \otimes \cW  \stackrel{f_i}{\rightarrow} \cF_i \rightarrow 0$$
by $\cO_X(k)$ and take the long exact sequence of cohomology, we
get an exact sequence
$$H \otimes W \stackrel{f_i}{\rightarrow} H^0(\cF_i(k)) \rightarrow H^1(\cK_i(k)).$$
The third term vanishes so that this gives
$$\dim f_i (H \otimes W) = \chi (\cF_i(k)).$$
Now Theorem \ref{thm:subspace}' can be applied
to the configuration
$$\{0 \rightarrow V \otimes W \stackrel{f_i}{\rightarrow} H^0(\cE_i(k)) \rightarrow 0 \}$$
to conclude the proof.
\endproof

\subsection{Moduli of Semistable Configuration of Coherent Sheaves}

Let ${\bfmit M}_P$ be the moduli space of semistable coherent sheaves over $X$ with
the Hilbert polynomial $P$.

Fix a set of positive numbers  ${\bf \omega} = \{\omega_1, \ldots, \omega_m\}$
and (Hilbert) polynomials ${\bf P} =\{P_1, \ldots, P_m\}$. Let
${\bfmit M}_{\bf P, \omega}$ be the moduli space of semistable configurations of coherent sheaves over $X$ with
the Hilbert polynomial $P_i$ and with respect to the weight ${\bf \omega} = \{\omega_1, \ldots, \omega_m\}$.

\begin{prop} Fix an integer $1 \le i \le m$. For sufficient large $\omega_i$ (relative
to other $\omega_j$), we have
\begin{enumerate}
\item If a configuration $\{V \otimes  \cW \rightarrow \cE_j \rightarrow 0\}$
of coherent sheaves is ${\bf \omega}$-semistable, then its $i$-th component
$V \otimes  \cW \rightarrow \cE_i \rightarrow 0$ is a semistable  sheaf;
\item If the $i$-th component
$V \otimes  \cW \rightarrow \cE_i \rightarrow 0$ of
 a configuration $\{V \otimes  \cW \rightarrow \cE_j \rightarrow 0\}$ is stable, then
the configuration $\{V \otimes  \cW \rightarrow \cE_j \rightarrow 0\}$ is ${\bf \omega}$-stable;
\item
In particular, there is a porjective morphism from ${\bfmit M}_{\bf P, \omega}$ to ${\bfmit M}_P$
$$\pi_i: {\bfmit M}_{\bf P, \omega} \rightarrow {\bfmit M}_P.$$
\end{enumerate}
\end{prop}
\proof
The proof is completely similar to that of Proposition \ref{projection}, thus is omitted.
\endproof

Recall that the stability of a coherent sheaf $\cE$ is defined as follows. $\cE$ is semi stable
(resp. stable) if for every proper subsheaf $\cF$ of $\cE$ we have that
$${\chi(\cF(k)) \over \hbox{rk}(\cF)} \le {\chi(\cE(k)) \over \hbox{rk}(\cE)}$$
(resp. $<$) for sufficiently large $k$ (e.g., \cite{Friedman},
\cite{Gieseker}, \cite{HL}, \cite{Simpson94}). It would be nice to
also have an intrinsic stability criterion (definition) for
configurations of coherent sheaves without using Grothendieck's
Grassmannian embeddings. Other directions of further research
include: to study the properties of the moduli (cf.,  e.g.,
\cite{GL} and \cite{Li}), and to study the dependence of the
moduli on the parameters (cf., e.g., \cite{FQ} and \cite{Qin},
among others).

\section{Balanced Configuration and Moment Map}

\subsection{Quot scheme and $\Hom(X,\Gr)$}
After tensoring coherent sheaves by $\cO_X(N)$ for large enough $N$,
we may assume that they are generated by global sections, hence regard
them as quotient sheaves of the trivial sheaf $V = {\Bbb C}^N \times X$
$$V  \stackrel{f_i}{\rightarrow} U_i \rightarrow 0.$$

We will focus on vector bundles only. This allows us to switch the viewpoint
and consider vector subbundles of ${\Bbb C}^N \times X$ instead of quotient bundles.
So, let
$$\cE_i \subset {\Bbb C}^N \times X$$
be a configuration of vector subbundles of rank $r_i$ over $X$ with
 the Hilbert polynomial $P_i$ ($1\le i \le m$).

Each $\cE_i$ corresponds to a map
$$g_i: X \rightarrow \Gr(r_i, {\Bbb C}^N )$$
where $g_i$ sends $x \in X$ to the fiber $(\cE_i)_x \subset {\Bbb C}^N$.
Conversely, every morphism
$$g: X \rightarrow \Gr(r, {\Bbb C}^N )$$
defines a vector subbundle by pulling back the universal bundle
$$\cE=\{ (v,x) \in {\Bbb C}^N \times X | v \in g(x)\}.$$
Let $$\Hom(X, \Gr(r_i, {\Bbb C}^N ); P_i)$$ be the set of morphisms that correspond to
vector subbundles of Hilbert polynomial $P_i$. Then we have an embedding
$$j: \Pi_{i=1}^m \Hom(X, \Gr(r_i, {\Bbb C}^N ); P_i) \rightarrow \Pi_{i=1}^m \Quot(V,P_i).$$
We will use the pull-back bundle $j^* L_{{\bf \omega}}$ as the linearization on
$$\Pi_{i=1}^m \Hom(X, \Gr(r_i, {\Bbb C}^N ); P_i)$$
where $L_{{\bf \omega}}$ is $L_{1, {\bf \omega}}$ as defined in \S 7.2.
Intrinsically, this bundle admits a description similar to $L_{m, {\bf \omega}}$.
Consider the diagram
\begin{equation*}
\begin{CD}
 \Hom (X, \Gr(r_i, {\Bbb C}^N )) \times X @>{ev}>> \Gr  \\
@V{\pi}VV  \\
 \Hom (X, \Gr(r_i, {\Bbb C}^N ))
\end{CD}
\end{equation*}
Let ${\cU}_i$ be the universal vector bundle over $\Gr$. Then
$$L_i=\Det(\pi_*(ev_*({\cU}_i \otimes  \cO_X(1)))$$
is very ample. For a weight set ${\bf \omega}$,
the  tensor product $\otimes_i L_i^{\omega_i}$ of these line bundles
on $\Pi_{i=1}^m \Hom(X, \Gr(r_i, {\Bbb C}^N ); P_i)$ is $j^* L_{{\bf \omega}}$.

By Lemma \ref{inclusionstable}, a configuration $\{\cE_i\}$ of vector subbundles of
${\Bbb C}^N \times X$ is (semi) stable
with respect to $L_{\bf \omega}$ if and only
if the corresponding configuration of morphisms
$g_i: X  \to \Gr(r_i, {\Bbb C}^N )$ is (semi) stable with respect to $j^*L_{\bf \omega}$.

\subsection{Moment map for singular varieties}

Let $Z$ be any (possibly) singular variety acted upon by a compact group $K$.
Let $\Omega$ be a bilinear skew-symmetric form on the Zariski tangent space $TZ$
which restricts to a symplectic form on $Z^0$, the smooth locus of $Z$. A continuous equivariant map
$$\Phi: Z \rightarrow {\bfmit k}^*$$
is called a moment map if the resrtiction
$$\Phi_{Z^0}: Z^0 \rightarrow {\bfmit k}^*$$
is a moment map (in the usual sense) for the action of $K$ on $Z^0$.
That is, at a smooth point of $Z$, we have
$$ d \langle \Phi, a \rangle = i_{\xi_a} \omega$$
for every $a \in {\bfmit k}$ where $\xi_a$ is the vector field generated by $a$.
By  continuity, the  moment map $\Phi: Z \rightarrow {\bfmit k}^*$
is uniquely determined by the
moment map $\Phi_{Z^0}: Z^0 \rightarrow {\bfmit k}^*.$
Note also that the moment map, when exists, is unique if the group $G$ is semisimple.
For linear actions on projective varieties, a moment map always exists.

If in addition, $Z$ can be equivariantly embedded in a smooth ambient variety $\widetilde{Z}$,
then the restriction of a moment map
$$\widetilde{\Phi}: \widetilde{Z} \rightarrow  {\bfmit k}^*$$
to $Z$ will be a moment map for the $K$-action on $Z$.
This situation is  the case we will be interested. That is, we will consider the
equivariant embeddings of the Quot schemes in the Grassmannians.
Lemma \ref{inclusionstable} will allow us to apply some results in the smooth case
to the singular case.

\subsection{Moment map for $\Pi_{i=1}^m \Hom(X, \Gr(r_i, {\Bbb C}^N ); P_i)$}

Now consider the space $$\Pi_{i=1}^m \Hom(X, \Gr(r_i, {\Bbb C}^N
); P_i).$$ $\SL(N)$ acts on it diagonally by moving the images. We
assume that $\Pi_{i=1}^m \Hom(X, \Gr(r_i, {\Bbb C}^N ); P_i)$ is
generically smooth (hence every component $\Hom(X, \Gr(r_i, {\Bbb
C}^N ); P_i)$ is also generically smooth). The line bundle
$j^*L_{{\bf \omega}}$ induces a symplectic form $\Omega$ on the
smooth locus of $\Pi_{i=1}^m \Hom(X, \Gr(r_i, {\Bbb C}^N ); P_i)$
as follows. At any given  point $f: X\hookrightarrow \Gr(r_i,
{\Bbb C}^N )$, the tangent space $T_f \Hom(X, \Gr(r_i, {\Bbb C}^N
))$
 is $H^0(X, f^*T\Gr(r_i, {\Bbb C}^N ))$.
We can define a skew-symmetric bilinear form $\Omega_i$ on $H^0(X, f^*T\Gr(r_i, {\Bbb C}^N ))$
by setting
$$(\Omega_i)_f (u,v) = \int_X f^*(\omega_i)_{FS} (u, v) dV$$
where $u, v \in H^0(X, f^*T \Gr(r_i, {\Bbb C}^N ))$ and
$(\omega_i)_{FS}$ is the symplectic form induced from the
Fubini-Study K\"ahler form on $\Gr(r_i, {\Bbb C}^N )$. The form
$\Omega_i$ restricts to a symplectic form on the smooth locus of
$\Hom(X, \Gr(r_i, {\Bbb C}^N ); P_i)$. Then the form on
$\Pi_{i=1}^m \Hom(X, \Gr(r_i, {\Bbb C}^N ); P_i)$ is
$$\Omega = \sum_{i=1}^m \omega_i \Omega_i.$$

Let $\Vol(X)$ be the volume of $X$ and $I$ denote the identity matrix in $\su(N)$.
Then we have
\begin{prop}
\label{momentmap}
Under the symplectic form $\Omega$, the moment map $\Phi$ of the action of $\SU(N)$ on
$\Pi_{i=1}^m \Hom(X, \Gr(r_i, {\Bbb C}^N ); P_i)$ is given by
$$\Phi (\{g_i\}) = \sum_{i=1}^m \omega_i \int_X A_i(x)A_i^*(x) d V - {\wp}_{\bf \omega}(\{g_i\})  \Vol(X) I$$
where $\{g_i\} \in \Pi_{i=1}^m \Hom(X, \Gr(r_i, {\Bbb C}^N ); P_i)$,
$A_i(x)$ is a matrix  representation of $g_i (x) \subset {\Bbb C}^N$ whose
columns is an orthonormal basis for $g_i(x)$ ($1 \le i \le m$), and
${\wp}_{\bf \omega}(\{g_i\}) = \sum_i \omega_i {r_i \over N}$.
\end{prop}

\proof
One first checks that for any given $i$ the moment map $\Phi_i$ of the action of $\SU(N)$ on
$\Hom(X, \Gr(r_i, {\Bbb C}^N ); P_i)$ is the integration over $X$
 of the moment map $\phi_i$  of  the  action of $\SU(N)$ on
the Grassmannian $\Gr(r_i, {\Bbb C}^N)$.
For any $a \in \su (N)$, it generates a vector field $\xi_a$  on $\Gr(r_i, {\Bbb C}^N)$.
At any smooth point $f \in \Hom(X, \Gr(r, {\Bbb C}^N ); P)$, we have
$$i_{\xi_a} (\Omega_i)_f = \int_X f^* i_{\xi_a} (\omega_i)_{FS} d V$$
$$=\int_X f^* \langle d \phi_i, a \rangle \wedge d V$$
$$= \pi_*(ev^*(d \langle  \phi_i, a \rangle) \wedge dV )$$
$$= d \langle \pi_* (ev^* \phi_i dV), a \rangle$$
$$= d \langle \int_X \phi_i dV, a \rangle.$$
This implies that $\Phi_i = \int_X \phi_i d V$.

Therefore the moment map  $\Phi$ of the action of $\SU(N)$ on
$$\Pi_{i=1}^m \Hom(X, \Gr(r_i, {\Bbb C}^N ); P_i)$$ is the same as the integration over $X$
 of the moment map $\Phi_0$ of  the diagoanl action of $\SU(N)$ on the product of
the Grassmannians.
Since
$$\Phi_0 (\{g_i (x)\})= \sum_i \omega_i (A_i(x)A_i^*(x) - {r_i \over N} I),$$
we have
$$\Phi (\{g_i \}) = \int_X \Phi_0 d \Vol=
\sum_i \int_X \omega_i (A_i(x)A_i^*(x) - {r_i \over N} I) d\Vol.$$
That is $$\Phi (\{g_i\}) = \sum_i \omega_i \int_X A_i(x)A_i^*(x) d V - {\wp}_{\bf \omega}(\{g_i\}) \Vol (X) I. $$
\endproof


\subsection{Balanced Configuration and Stability}

\begin{defn}
Let $\{g_i: X \rightarrow \Gr(r_i, {\Bbb C}^N)\}$ be a configuration of
morphism into the Grassmannians.
We say that the configuration $\{g_i\}$ is balanced if
$$\sum_{i=1}^m \omega_i \int_X A_i(x)A_i^*(x) d V = {\wp}_{\bf \omega}(\{g_i\})  \Vol(X) I.$$
We say $\{g_i\}$ can be (uniquely) balanced if there is a (unique) element $u \in \SU(N) \backslash SL(N)$ such that
$\{u \cdot g_i\}$ is balanced.
\end{defn}

The following theorem follows from Proposition \ref{momentmap} and Lemma \ref{inclusionstable}.

\begin{thm}
\label{stable=balanced}
$\{g_i\}$ is stable if and only if $\{g_i\}$ can be (uniquely) balanced and its stabilizer
group is finite.
\end{thm}

\begin{defn}
Let $\{\cE_i\}$ be a configuration of vector subbundles in
$\Pi_{i=1}^m \Quot(V,P_i)$ where $V$ is the trivial vector bundle
${\Bbb C}^N \times X$. We say  the system $\{\cE_i\}$ is balanced
if
$$\sum_{i=1}^m \omega_i \int_X A_i(x)A_i^*(x) d V = {\wp}_{\bf \omega}(\{\cE_i\}) \Vol(X) I$$
where $A_i(x)$ is a matrix  representation of $(\cE_i)_x \subset {\Bbb C}^N$ whose
columns is an orthonormal basis for $(\cE_i)_x$ ($1 \le i \le m$), and
${\wp}_{\bf \omega}(\{\cE_i\}) = \sum_i \omega_i {r_i \over N}$. .
We say $\{\cE_i\}$ can be (uniquely) balanced if there is a (unique) element $u \in \SU(N) \backslash \SL(N)$ such that
$\{u \cdot \cE_i\}$ is balanced.
\end{defn}

As a consequence of Theorem \ref{stable=balanced}, we obtain

\begin{thm}
\label{stable=balancedbundles}
Let $\{\cE_i\}$ be a configuration of vector subbundles in $\Pi_{i=1}^m \Quot(V,P_i)$. Then
$\{\cE_i\}$ is stable if and only if $\{\cE_i\}$ can be (uniquely) balanced and its stabilizer
group is finite.
\end{thm}

When $m >1$, the condition that ``the stabilizer group of the configuration
$\{\cE_i\}$ is finite'' is a quite weak condition. For example, it will be the case
when $\cap_i \Stab (\cE_i)$ is finite
where $\Stab (\cE_i)$ is the stabilizer group of $\cE_i$ ($1 \le i \le m$).

Finally, consider the case when $m=1$.
Let $\cE$ be a vector subbundle in ${\Bbb C}^N \times X$. Then we obtain
a result of Wang (\cite{Wang}) and Phong-Sturm (\cite{PhS})

\begin{thm}
$\cE$ is Gieseker-Simpson
 stable if and only if it can be (uniquely) balanced and its automorphism group is finite.
\end{thm}

\bibliographystyle{amsplain}

\begin{thebibliography}{10}
\bibitem{Borcea} C. Borcea,
{\em  Association for flag configurations.}  In: ``Commutative
Algebra, Singularities and Computer Algebra'', NATO Science Series
II, Vol. 115, Kluwer Ac. Publ., 2003.
\bibitem{Dolgachev} I. Dolgachev,
{\em Introduction to Geometric Invariant Theory Quotients},
Lecture Notes Series, No. 25, Seoul National Univ. 1994.
\bibitem{Dolgachev2003} I. Dolgachev,
{\em Lectures on invariant theory.} London Mathematical Society
Lecture Note Series, {\bf 296}. Cambridge University Press,
Cambridge, 2003.
\bibitem{DH} I. Dolgachev and Y. Hu:
{\em Variation of Geometric Invariant Theory}, with an appendix by
Nicolas Ressayr. Publ. Math. I.H.E.S. {\bf 78} (1998), 1 -- 56.
\bibitem{Donaldson} S. Donaldson,
{\em Scalar curvature and projective embeddings. I,}
 J. Differential  Geom. {\bf 59} (2001), no. 3, 479--522.
 \bibitem{EP} D. Eisenbud and S. Popescu,
 {\em The projective geometry of the Gale transform.}
  J. Algebra  {\bf 230}  (2000),  no. 1, 127--173.
\bibitem{Faltings} G. Faltings and G. W\"ustholz,
{\em Diophantine approximations on projective spaces,}
Invent. Math {\bf 116} (1994), 109--138.
\bibitem{Foth} P. Foth and G. Lozano,
{\em The Geometry of Polygons in ${\Bbb R}^5$ and Quaternions,}
         math.DG/0202162
\bibitem{Friedman} R. Friedman,
{\em Algebraic surfaces and holomorphic vector bundles.}
Universitext. Springer-Verlag, New York, 1998.
\bibitem{FQ} R. Friedman and Z. Qin,
{\em Flips of moduli spaces and transition formulas for Donaldson
polynomial invariants of rational surfaces.} Comm. Anal. Geom.  3
(1995),  no. 1-2, 11--83.
\bibitem{GM} I. Gelfand and R. MacPherson,
{\em Geometry in Grassmannians and a generalization of the
dilogarithm. }  Adv. in Math.  44  (1982), no. 3, 279--312.
\bibitem{Gieseker} D. Gieseker,
{\em On the moduli of vector bundles on an algebraic surface.}
Ann. of Math. (2)  106  (1977), no. 1, 45--60.
\bibitem{GL} D. Gieseker and J. Li,
{\em  Irreducibility of moduli of rank-$2$ vector bundles on
algebraic surfaces.}  J. Differential Geom.  40  (1994),  no. 1,
23--104.
\bibitem{Hu96} Y. Hu,
{\em Relative Geometric Invariant Theory and Universal Moduli Spaces,}
International Journal of Mathematics  Vol. 7 No. 2 (1996), 151 -- 181.
\bibitem{Hu2003} Y. Hu,
{\em Topological Aspects of Chow Quotients.} math.AG/0308027.
\bibitem{HL} D. Huybrechts and M. Lehn,
{\em The geometry of moduli spaces of sheaves.}  Aspects of
Mathematics, E31. Friedr. Vieweg \& Sohn, Braunschweig, 1997.
\bibitem{Kapranov95} M. Kapranov,
{\em Chow quotients of Grassmannian,} I.M. Gelfand Seminar Collection, AMS.
\bibitem{Kly} A. Klyachko,
{\em Stable bundles, representation theory and Hermitian operators,}
Selecta Math. {\bf 4} (1998), 419-445.
\bibitem{Li} J. Li,
{\em Algebraic geometric interpretation of Donaldson's polynomial
invariants.}  J. Differential Geom.  37  (1993),  no. 2, 417--466.
\bibitem{Luo} H.-Z. Luo,
{\em Geometric criterion for Gieseker-Mumford stability of polarized manifolds,}
Journal of Diff. Geom. {\bf 49} (1998), no. 3, 577--599.
\bibitem{GIT} D. Mumford, J. Fogarty, and F. Kirwan,
Geometric Invariant Theory,
Springer-Verlag, Berlin, New York,  1994.
\bibitem{PhS} D.H. Phong and J. Sturm,
{\em Stability, Energy Functions, and K\"ahler-Einstein Metrics,}
arXiv:math.DG/0203254.
\bibitem{Qin} Z. Qin,
{\em Equivalence classes of polarizations and moduli spaces of
sheaves.}  J. Differential Geom.  37  (1993),  no. 2, 397--415.
\bibitem{Simpson94} C. Simpson,
{\em Moduli spaces of representations of the fundamental group of
a smooth projective variety, I and II},  Inst. Hautes Études Sci.
Publ. Math.  No. 79 (1994), 47--129.  Inst. Hautes Études Sci.
Publ. Math.  No. 80 (1994), 5--79.
\bibitem{Tian} G. Tian,
{\em On a set of polarized K\"ahler metrics on algebraic manifolds,}
Jour. Diff. Geom. {\bf 32} (1990), 99-130.
 \bibitem{Tian91} G. Tian,
{\em Kähler-Einstein metrics on algebraic manifolds.}
 Proceedings of the International Congress of
  Mathematicians, Vol. I, II (Kyoto, 1990), 587--598, Math. Soc. Japan, Tokyo, 1991.
\bibitem{Tian97} G. Tian,
{\em Kähler-Einstein metrics with positive scalar curvature.}
Invent. Math. {\bf 130} (1997), 1--39.
\bibitem{Totaro} B. Totaro,
{\em Tensor products of semistables are semistable,}
Geometry and Analysis on complex manifolds, World Sci. Publ. 1994, pp. 242-250.
\bibitem{Viehweg95} E. Vieweg,
{\em Quasi-projective moduli of polarized manifolds,}
Springer-Verlag, 1995.
\bibitem{Wang} X. W. Wang,
{\em Balanced point, stability and vector bundles over projective manifolds,}
preprint.
\bibitem{Yau} S.-T. Yau,
{\em Open problems in geometry,}
Proc. Symp. Pure Math. {\bf 54} (1993) 1-28.
\bibitem{zhang} S. Zhang,
{\em Heights and reductions of semi-stable varieties,}
Compositio Mathematica {\bf 104} (1996) 77-105.
\end{thebibliography}
\makeatletter \renewcommand{\@biblabel}[1]{\hfill#1.}\makeatother

\end{document}